# CONSERVATIVE STOCHASTIC CAHN–HILLIARD EQUATION WITH REFLECTION

By Arnaud Debussche and Lorenzo Zambotti

*Antenne De Bretagne and Université Paris 6*

We consider a stochastic partial differential equation with reflection at 0 and with the constraint of conservation of the space average. The equation is driven by the derivative in space of a space–time white noise and contains a double Laplacian in the drift. Due to the lack of the maximum principle for the double Laplacian, the standard techniques based on the penalization method do not yield existence of a solution. We propose a method based on infinite dimensional integration by parts formulae, obtaining existence and uniqueness of a strong solution for all continuous nonnegative initial conditions and detailed information on the associated invariant measure and Dirichlet form.

**1. Introduction.** Consider the following stochastic Cahn–Hilliard equation on $[0,1]$ with homogeneous Neumann boundary condition and reflection at $u = 0$:

$$(1.1) \quad \begin{cases} \dfrac{\partial u}{\partial t} = -\dfrac{1}{2}\dfrac{\partial^2}{\partial \theta^2}\left(\dfrac{\partial^2 u}{\partial \theta^2} + \eta\right) + \dfrac{\partial}{\partial \theta}\dot{W}, \\ \dfrac{\partial u}{\partial \theta}(t,0) = \dfrac{\partial u}{\partial \theta}(t,1) = \dfrac{\partial^3 u}{\partial \theta^3}(t,0) = \dfrac{\partial^3 u}{\partial \theta^3}(t,1) = 0, \\ u(0,\theta) = u_0(\theta), \qquad \theta \in [0,1], \end{cases}$$

where $\dot{W}$ is a space–time white noise on $[0,+\infty) \times [0,1]$, $u$ is a continuous function of $(t,\theta) \in [0,+\infty) \times [0,1]$, $\eta$ is a locally finite positive measure on $(0,+\infty) \times [0,1]$, subject to the constraint

$$(1.2) \quad u \geq 0, \qquad \int_{(0,+\infty)\times[0,1]} u\, d\eta = 0.$$

Stochastic partial differential equations with reflection can model the evolution of random interfaces near a hard wall; see [17] and [26]. However in these papers the equation of interest contains a second-order operator rather than a fourth-order one. More precisely, the studied equation is

$$
\text{(1.3)} \quad \begin{cases} \dfrac{\partial v}{\partial t} = \dfrac{1}{2}\dfrac{\partial^2 v}{\partial \theta^2} + \zeta + \dot{W}, \\ v_t(0) = v_t(1) = 0, \qquad t \geq 0, \\ v \geq 0, \qquad d\zeta \geq 0, \qquad \int v\, d\zeta = 0, \end{cases}
$$

where $v$ is a continuous function and $\zeta$ a positive measure on $\mathcal{O}$. In [17] it is proven that the fluctuations of a $\nabla \phi$ interface model near a hard wall converge in law to the stationary solution of (1.3). However, if one introduces the constraint of conservation of the area between the interface and the wall, then the Cahn–Hilliard equation (1.1) is expected; see [23] and [17]. For other results on fluctuations of random interfaces, see [18].

Equation (1.3) has been introduced by Nualart and Pardoux in [21] and studied in detail recently: besides existence and uniqueness of solutions, a rather detailed study of the contact set $\{(t,\theta): v_t(\theta) = 0\}$ and of the reflection measure $\zeta$ has been obtained; see [7, 24, 25].

A crucial tool in these papers, including [17], is the following monotonicity property: if we define $v^\varepsilon$ as the unique solution of

$$\frac{\partial v^\varepsilon}{\partial t} = \frac{1}{2}\frac{\partial^2 v^\varepsilon}{\partial \theta^2} + \frac{1}{\varepsilon}(v^\varepsilon)^- + \dot{W}, \qquad v_0^\varepsilon = u_0,$$

with Dirichlet boundary condition at 0 and 1, then

$$\text{(1.4)} \qquad\qquad 0 < \varepsilon \leq \varepsilon' \quad \Longrightarrow \quad v^\varepsilon \geq v^{\varepsilon'},$$

that is, $\varepsilon \mapsto v^\varepsilon$ is a monotone function. This fundamental property stems from the maximum principle of the Laplace operator and is used, for example, in [21] to prove existence of a solution of (1.3).

The classical penalization method is used also in this paper: indeed, we introduce an approximation of (1.1) by means of the following stochastic partial differential equation with a Lipschitz nonlinear term:

$$
\text{(1.5)} \quad \begin{cases} \dfrac{\partial u^\varepsilon}{\partial t} = -\dfrac{1}{2}\dfrac{\partial^2}{\partial \theta^2}\bigg(\dfrac{\partial^2 u^\varepsilon}{\partial \theta^2} + \dfrac{1}{\varepsilon}f(u^\varepsilon)\bigg) + \dfrac{\partial}{\partial \theta}\dot{W}, \\ \dfrac{\partial u^\varepsilon}{\partial \theta}(t,0) = \dfrac{\partial u^\varepsilon}{\partial \theta}(t,1) = \dfrac{\partial^3 u^\varepsilon}{\partial \theta^3}(t,0) = \dfrac{\partial^3 u^\varepsilon}{\partial \theta^3}(t,1) = 0, \\ u^\varepsilon(0,\theta) = u_0(\theta), \end{cases}
$$

where $f : \mathbb{R} \mapsto \mathbb{R}$ is monotone nonincreasing and Lipschitz-continuous with $f(u) = 0$ for $u \geq 0$ and $f(u) > 0$ for $u < 0$, for instance,

$$f(u) = -(u \wedge 0) = (u)^-, \qquad u \in \mathbb{R}.$$

However, the monotonicity property (1.4) which holds for (1.3) fails for (1.1), that is, in general there exist no $\varepsilon \neq \varepsilon'$ such that $u^\varepsilon \geq u^{\varepsilon'}$. Indeed, integrating in $d\theta$ over $[0, 1]$ equation (1.5) we obtain that

$$(1.6) \qquad \int_0^1 u^\varepsilon(t, \theta) \, d\theta = \int_0^1 u_0(\theta) \, d\theta = \int_0^1 u^{\varepsilon'}(t, \theta) \, d\theta \qquad \forall \, t \geq 0,$$

and if $u^\varepsilon \geq u^{\varepsilon'}$ then necessarily $u^\varepsilon \equiv u^{\varepsilon'}$, which is impossible if $\varepsilon \neq \varepsilon'$. Therefore $\varepsilon \mapsto u^\varepsilon$ is not monotone and the techniques used for (1.3) break down. In fact even well-posedness of (1.1) appears to be a new result.

In this paper, we propose an approach to (1.1) which yields well-posedness but also a description of the invariant measures and an integration by parts formula on such measures which gives a better insight into the reflection measure $\eta$.

Although we do not have monotonicity of $\varepsilon \mapsto u^\varepsilon$, we do want to prove that $u^\varepsilon$ converges in a suitable sense to the unique solution $u$ of (1.1). We first notice that (1.5) is a gradient system in $H^{-1}(0, 1)$ with a convex potential (see Section 3 below). This fact yields in particular the crucial strong Feller property uniformly in $\varepsilon > 0$, that is, the equicontinuity of the transition semigroup of the Markov process $(u_t^\varepsilon, t \geq 0)$ in the topology of $H^{-1}(0, 1)$.

Analogously to (1.6), the solutions of (1.1) have a constant space average $\int_0^1 u_t \, d\theta = \int_0^1 u_0 \, d\theta =: c$, $t \geq 0$. If $c > 0$ is fixed, we prove tightness of the stationary solutions of (1.5) thanks to the Lyons–Zheng decomposition and by pathwise uniqueness of solutions of (1.1), we have a unique limit. Thus, we obtain a unique stationary solution of (1.1) and the invariant measure $\nu_c$ can be explicitly described. Then, one obtains existence of solutions of (1.1) for $\nu_c$-a.e. initial condition $u_0$. Moreover, one can prove that one can extend this process by density to a continuous process in the space of distributions $H^{-1}(0, 1)$ for any initial condition $u_0$.

However, we have no way now to give a meaning to the equation for this process, since the contact condition (1.2) requires a priori continuity of the solution, while we know only that $u_t(\cdot) \in H^{-1}(0, 1)$. We solve this problem proving a much stronger statement than convergence of stationary measures of $u^\varepsilon$: we prove that the transition semigroup of the Markov process $(u_t^\varepsilon, t \geq 0)$ converges in a suitable sense to the Markov transition semigroup of $(u_t, t \geq 0)$, which therefore inherits the strong Feller property. The latter property allows in turn to prove that the process constructed above is continuous on $(0, \infty) \times [0, 1]$ for any initial condition.

We remark that the techniques of this paper can be applied to prove existence of solutions of a Cahn–Hilliard equation

$$\begin{cases} \dfrac{\partial u}{\partial t} = -\dfrac{1}{2} \dfrac{\partial^2}{\partial \theta^2} \left( \dfrac{\partial^2 u}{\partial \theta^2} + g(u) \right) + \dfrac{\partial}{\partial \theta} \dot{W}, \\ \dfrac{\partial u}{\partial \theta}(t, 0) = \dfrac{\partial u}{\partial \theta}(t, 1) = \dfrac{\partial^3 u}{\partial \theta^3}(t, 0) = \dfrac{\partial^3 u}{\partial \theta^3}(t, 1) = 0, \\ u(0, \theta) = u_0(\theta), \end{cases}$$

where $g:\mathbb{R} \mapsto \mathbb{R}$ is an arbitrary continuous function such that $u \mapsto g(u) + bu$ is monotone nonincreasing for some $b \in \mathbb{R}$. For the reason discussed above, existence of a solution is not trivial: indeed, existence in $H^{-1}(0,1)$ is not enough to give a meaning to the nonlinear term $g(u)$. However a suitable uniform strong Feller property and then existence of continuous solutions can be proven with the techniques of this paper.

We notice that the results of [7, 24] on the random contact set $\{(t,\theta) : v_t(\theta) = 0\}$ of equation (1.3) could have natural and interesting generalizations to equation (1.1) but look quite challenging. In this direction, we study the Dirichlet form associated with the solution of (1.1) and a related integration by parts formula for the invariant measure $\nu_c$, analogous to that computed in [24] and applied to the study of (1.3). For a recent result on infinite dimensional integration by parts formulae with boundary terms, see [16].

The stochastic Cahn–Hilliard equation is also a model arising in the study of spinodal decomposition. It has been the object of several mathematical studies dealing with well-posedness [8, 9, 13], properties of the solutions [5] or qualitative behavior [1, 2]. In these articles, the equation studied contains an extra nonlinear term as the one above. It is never required that the solution be positive and no reflexion term is necessary. The stochastic Cahn–Hilliard describes the evolution of a concentration and some authors consider a nonlinear term which is singular at $-1$ and $1$, a typical example being $g(u) = \ln(\frac{1+u}{1-u})$. In the deterministic case, when such a model is considered, the solution takes values in $[-1, 1]$ [3, 11]. When an additive noise term is considered, the singularity is not strong enough to prevent the solution from reaching the values $-1$ or $1$ and a reflexion term naturally arises. Our techniques extend to this type of equation, with additional difficulties due to the fact that two constraints are considered.

The plan of the paper is the following: we discuss in Section 2 the linear Cahn–Hilliard equation, in Section 3 the basic properties of (1.5), in Section 4 pathwise uniqueness, in Section 5 existence of stationary solutions of (1.1) and in Section 6 existence of strong solutions of (1.1). In Section 7 we state the integration by parts for the invariant measures of (1.1), whose proof is postponed to Appendix A; in Section 8 we study a related Dirichlet form.

1.1. *Notation.* We denote by $\langle \cdot, \cdot \rangle_L$ the canonical scalar product in $L^2(0,1)$:
$$\langle h, k \rangle_L := \int_0^1 h_\theta k_\theta \, d\theta.$$

We denote by $A$ the realization in $L^2(0,1)$ of $\partial_\theta^2$ with Neumann boundary condition at 0 and 1, that is,

(1.7) $\quad D(A) := \{h \in H^2(0,1) : h'(0) = h'(1) = 0\}, \qquad A := \dfrac{\partial^2}{\partial \theta^2}.$

Notice that $A$ is self-adjoint in $L^2(0,1)$, with the complete orthonormal system $(e_i)_i$ of eigenvectors in $L^2(0,1)$,

$$(1.8) \quad e_0(\theta) := 1, \qquad e_n(\theta) := \sqrt{2}\cos(n\pi\theta), \qquad Ae_n = -(n\pi)^2, \qquad n \in \mathbb{N}.$$

We introduce a notation for the *average* of $h \in L^2(0,1)$:

$$\overline{h} := \int_0^1 h = \langle h, e_0 \rangle_L.$$

Then we also set for all $c \in \mathbb{R}$

$$L_c^2 := \{h \in L^2(0,1) : \overline{h} = c\}.$$

Notice that $(-A)^{-1} : L_0^2 \mapsto L_0^2$ is well defined. More precisely, a direct computation shows that for all $h \in L^2(0,1)$,

$$Qh(\theta) := (-A)^{-1}[h - \overline{h}](\theta) = \int_0^1 q(\theta,\sigma) h_\sigma \, d\sigma,$$

where

$$(1.9) \qquad q(\theta,\sigma) := \theta \wedge \sigma + \frac{\theta^2 + \sigma^2}{2} - \theta - \sigma + \frac{1}{3}, \qquad \theta,\sigma \in [0,1].$$

We extend $Q$ to a one-to-one operator $\overline{Q} : L^2(0,1) \mapsto L^2(0,1)$ by

$$(1.10) \qquad \overline{Q}h := Q(h - \overline{h}) + \overline{h},$$

that is, for all $h \in L^2(0,1)$,

$$(1.11) \qquad \overline{Q}h(\theta) = \int_0^1 [q(\theta,\sigma) + 1] h_\sigma \, d\sigma.$$

Then we define $H$ as the completion of $L^2(0,1)$ with respect to the scalar product

$$(h,k)_H := \langle \overline{Q}h, k \rangle_L.$$

For all $c \in \mathbb{R}$ we also set

$$H_c := \{h \in H : (h, e_0)_H = c\}.$$

We remark that $H$ is naturally interpreted as a space of distributions, in particular as the dual space of $H^1(0,1)$. Finally, we denote by $\Pi$ the symmetric projector of $H$ onto $H_0$, that is,

$$(1.12) \qquad \Pi : H \mapsto H_0, \qquad \Pi h := h - \overline{h}, \qquad \overline{h} := (h, e_0)_H.$$

Notice that $\Pi$ is also a symmetric projector of $L^2(0,1)$ onto $L_0^2$.

1.2. *Weak and strong solutions of* (1.1). We state now the precise meaning of a solution to (1.1).

DEFINITION 1.1. Let $u_0 \in C([0,1])$, $u_0 \geq 0$, $\int_0^1 u_0 > 0$. We say that $(u, \eta, W)$, defined on a filtered complete probability space $(\Omega, \mathbb{P}, \mathcal{F}, \mathcal{F}_t)$, is a *weak* solution to (1.1) on $[0,T]$ if:

1. a.s. $u \in C((0,T] \times [0,1])$, $u \geq 0$ and $u \in C([0,T]; H)$;
2. a.s. $\eta$ is a positive measure on $(0,T] \times [0,1]$, such that $\eta([\delta, T] \times [0,1]) < \infty$ for all $\delta \in (0,T]$;
3. $(W(t, \theta))$ is a Brownian sheet, that is, a centered Gaussian process such that
$$\mathbb{E}[W(t,\theta) W(t',\theta')] = t \wedge t' \cdot \theta \wedge \theta', \qquad t, t' \geq 0, \theta, \theta' \in [0,1];$$
4. $u_0$ and $W$ are independent and the process $t \mapsto (u_t(\theta), W(t, \theta))$ is $(\mathcal{F}_t)$-adapted for all $\theta \in [0,1]$,
5. for all $h \in D(A^2)$ and for all $0 < \delta \leq t \leq T$,

(1.13)
$$\langle u_t, h \rangle_L = \langle u_\delta, h \rangle_L - \tfrac{1}{2} \int_\delta^t \langle u_s, A^2 h \rangle_L \, ds$$
$$- \tfrac{1}{2} \int_\delta^t \int_0^1 Ah_\theta \, \eta(ds, d\theta) - \int_\delta^t \int_0^1 h'_\theta \, W(ds, d\theta);$$

6. a.s. the contact property holds: $\operatorname{supp}(\eta) \subset \{(t,\theta) : u_t(\theta) = 0\}$, that is,
$$\int_{(0,T] \times [0,1]} u \, d\eta = 0.$$

Finally, a weak solution $(u, \eta, W)$ is a strong solution if the process $t \mapsto u_t(\cdot)$ is adapted to the filtration $t \mapsto \sigma(W(s, \cdot), u_0(\cdot) : s \in [0,t])$.

In Theorem 6.2 we shall prove that for all $u_0 \in C([0,1])$ with $u_0 \geq 0$ and $\int_0^1 u_0 > 0$ there exists a unique strong solution of (1.1). We shall also study the ergodic properties of the solutions, the associated transition semigroup and Dirichlet form.

1.3. *Function spaces.* Notice that for all $c \in \mathbb{R}$, $H_c = ce_0 + H_0$ is a closed affine subspace of $H$ isomorphic to the Hilbert space $H_0$. If $J$ is a closed affine subspace of a Hilbert space, we denote by $C_b(J)$, respectively $C_b^1(J)$, the space of all bounded uniformly continuous functions on $J$, respectively bounded and uniformly continuous together with the first Fréchet derivative. We also denote by $\operatorname{Lip}(J)$ the set of all $\varphi \in C_b(J)$ such that
$$[\varphi]_{\operatorname{Lip}(J)} := \sup_{h \neq k} \frac{|\varphi(h) - \varphi(k)|}{\|h - k\|_J} < \infty.$$

For all $\varphi \in C_b^1(H)$ we denote by $\partial_h \varphi$ the directional derivative of $\varphi$ along $h \in H$:

$$\partial_h \varphi(x) := \lim_{t \to 0} \frac{1}{t}(\varphi(x+th) - \varphi(x)), \qquad x \in H.$$

Notice that we have natural inclusions $C_b(H) \subset C_b(L^2(0,1))$, respectively, $C_b^1(H) \subset C_b^1(L^2(0,1))$. In particular, by the definition of gradients we have for $\varphi \in C_b^1(H) \subset C_b^1(L^2(0,1))$

$$\nabla \varphi : H \mapsto L^2(0,1), \qquad \langle \nabla \varphi(x), h \rangle_L = \partial_h \varphi(x),$$
$$\nabla_H \varphi : H \mapsto H, \qquad (\nabla_H \varphi(x), h)_H = \partial_h \varphi(x), \qquad x, h \in H.$$

Finally, we define $\mathrm{Exp}_A(H) \subset C_b(H)$ as the linear span of $\{\cos((h,\cdot)_H), \sin((h,\cdot)_H) : h \in D(A^2)\}$. In particular, by the definition of the scalar product in $H$,

(1.14) $\qquad \nabla_H \varphi = (e_0 \otimes e_0 - A) \nabla \varphi \qquad \forall \varphi \in \mathrm{Exp}_A(H).$

**2. The linear equation.** We start with the linear Cahn–Hilliard equation:

(2.1)
$$\begin{cases} \dfrac{\partial z}{\partial t} = -\dfrac{1}{2}\dfrac{\partial^4 z}{\partial \theta^4} + \dfrac{\partial}{\partial \theta}\dot W, \\ \dfrac{\partial z}{\partial \theta}(t,0) = \dfrac{\partial z}{\partial \theta}(t,1) = \dfrac{\partial^3 z}{\partial \theta^3}(t,0) = \dfrac{\partial^3 z}{\partial \theta^3}(t,1) = 0, \\ z_0(\theta) = 0. \end{cases}$$

The unique solution has an explicit representation in Fourier series

$$z_t(\theta) = \sum_{n=1}^{\infty} n\pi e_n(\theta) \int_0^t e^{-(t-s)n^4\pi^4/2}\, dw_s^n, \qquad t \geq 0, \theta \in [0,1],$$

where for $n \in \mathbb{N}$

$$\varepsilon_n(\theta) := \sqrt{2}\sin(n\pi\theta) = -\frac{1}{n\pi}e_n'(\theta), \qquad w_t^n := \int_0^1 \varepsilon_n(\theta) W(t, d\theta),$$

and $(w^n)_{n \in \mathbb{N}}$ is an independent sequence of standard Brownian motions. Clearly, $z$ is a Gaussian process. It is easy to prove that for all $T \geq 0$ there exists a constant $C > 0$ such that

$$\mathbb{E}[|z_t(\theta) - z_{t'}(\theta')|^2] \leq C(|t-t'|^{1/4} + |\theta - \theta'|), \qquad t, t' \in [0,T], \theta \in [0,1].$$

In particular, by Kolmogorov's criterion,

$$\text{a.s.} \quad z \in C^{1/8-\varepsilon, 1/2-\varepsilon}([0,T] \times [0,1]) \qquad \forall\, \varepsilon \in (0, 1/8).$$

More generally, we can introduce the Ornstein–Uhlenbeck process

(2.2) $\qquad \begin{cases} dZ_t = -\frac{1}{2}A^2 Z\, dt + B\, dW_t, \\ Z_0(x) = x \in L^2(0,1), \end{cases}$

where $W$ is a cylindrical white noise in $L^2(0,1)$ and

$$D(B) := H^1_0(0,1), \qquad B := \frac{d}{d\theta}, \qquad D(B^*) := H^1(0,1), \qquad B^* := -\frac{d}{d\theta},$$

and we notice that $BB^* = -A$. Then it is well known that $Z$ is equal to

$$Z_t(x) = e^{-tA^2/2}x + \int_0^t e^{-(t-s)A^2/2} B\, dW_s = e^{-tA^2/2}x + z_t(\cdot)$$

and that this process belongs to $C([0,\infty); L^2(0,1))$. Notice that

$$(2.3)\quad \langle Z_t(x), e_0\rangle_L = \langle x, e^{-tA^2/2}e_0\rangle_L + \int_0^t \langle B^* e^{-(t-s)A^2/2}e_0, dW_s\rangle_L = \langle x, e_0\rangle_L,$$

since $e^{-tA^2/2}e_0 = e_0$ and $B^* e^{-(t-s)A^2/2}e_0 = B^* e_0 = 0$. In particular, the average of $Z$ is constant. Now, the $L^2(0,1)$-valued r.v. $Z_t(x)$ has law

$$Z_t(x) \sim \mathcal{N}(e^{-tA^2/2}x, Q_t),$$

$$Q_t = \int_0^t e^{-sA^2/2} BB^* e^{-sA^2/2}\, ds = (-A)^{-1}(I - e^{-tA^2}).$$

If we let $t \to \infty$, the law of $Z_t(x)$ converges to the Gaussian measure on $L^2(0,1)$

$$(2.4)\qquad\qquad \mu_c := \mathcal{N}(ce_0, Q), \qquad c := \bar{x},$$

with covariance operator $Q$ and mean $ce_0$. Notice that the kernel of $Q$ is $\{te_0 : t \in \mathbb{R}\}$ and that $\mu_c$ is concentrated on the affine space $L^2_c$. Finally, we introduce the Gaussian measure on $L^2(0,1)$

$$(2.5)\qquad\qquad \mu := \mathcal{N}(0, \overline{Q});$$

recall (1.11). In this case, the kernel of $\overline{Q}$ in $L^2(0,1)$ is the null space, so the support of $\mu$ is the full space $L^2(0,1)$. The next result gives a description of $\mu$ and $\mu_c$ as laws of stochastic processes related to the Brownian motion.

LEMMA 2.1. *Let $(B_\theta)_{\theta \in [0,1]}$ a Brownian motion and $a \in \mathcal{N}(0,1)$, such that $\{B, a\}$ are independent. If we set*

$$Y_\theta := B_\theta - \overline{B} - a, \qquad \overline{B} := \int_0^1 B,$$

$$Y^c_\theta := B_\theta - \overline{B} + c, \qquad \theta \in [0,1],$$

*then the law of $Y$ is $\mu$ and the law of $Y^c$ is $\mu_c$.*

PROOF. Clearly $Y$ is a centered Gaussian process. A computation shows that its covariance function is, for all $\theta, \sigma \in [0,1]$,

$$\mathbf{E}[Y_\theta Y_\sigma] = \theta \wedge \sigma + \frac{\theta^2 + \sigma^2}{2} - \theta - \sigma + \frac{4}{3} = q(\theta, \sigma) + 1.$$

By (1.9), (1.11) and (2.5) we have that $\mu$ is the law of $Y$. Analogously, $Y^c$ is a Gaussian process with mean $c$ and covariance function $q$, which proves the second assertion. $\square$

Since $a$ and $B - \overline{B}$ are independent, and $\overline{Y} = -a$, then we obtain that $\mu_c$ is a regular conditional distribution of $\mu(dx)$ given $\{\bar{x} = c\}$ for all $c \in \mathbb{R}$, that is,

$$\mu_c(dx) = \mu(dx|\bar{x} = c) = \mu(dx|L_c^2).$$

Recall (1.12) and (1.14). Then we have the following result:

PROPOSITION 2.2. *Let $c \in \mathbb{R}$. The bilinear form*

$$\begin{aligned}\Lambda^c(\varphi, \psi) &:= \tfrac{1}{2} \int_H (\Pi \nabla_H \varphi, \nabla_H \psi)_H \, d\mu_c \\ &= \tfrac{1}{2} \int_H \langle -A \nabla \varphi, \nabla \psi \rangle_L \, d\mu_c \qquad \forall \, \varphi, \psi \in \mathrm{Exp}_A(H),\end{aligned}$$

*is closable in $L^2(\mu_c)$ and the process $(Z_t(x) : t \geq 0, x \in H_c)$ is associated with the resulting symmetric Dirichlet form $(\Lambda^c, D(\Lambda^c))$. Moreover,* $\mathrm{Lip}(H_c) \subset D(\Lambda^c)$ *and* $\Lambda^c(\varphi, \varphi) \leq [\varphi]^2_{\mathrm{Lip}(H_c)}$.

PROOF. The proof is standard, since the process $Z$ is Gaussian; see [10], Chapter 10.2. However we include some details since the interplay between the topologies of $H$ and $L^2(0,1)$ can produce some confusion. The starting point is the following integration by parts formula for $\mu$:

$$(2.6) \quad \int \partial_h \varphi \, d\mu = \int \langle \overline{h} - Ah, x \rangle_L \varphi(x) \mu(dx) = \int (\overline{h} + A^2 h, x)_H \varphi(x) \mu(dx)$$

for all $\varphi \in C_b^1(H)$ and $h \in D(A^2)$. Recall that $\bar{x}$ and $\Pi x$ are independent under $\mu(dx)$. Then (2.6) implies

$$(2.7) \quad \int \partial_{\Pi h} \varphi \, d\mu_c = \int \langle -Ah, x \rangle_L \varphi(x) \mu_c(dx) = \int (A^2 h, x)_H \varphi(x) \mu_c(dx).$$

Let now $\varphi(x) := \exp(i(x,h)_H)$ and $\psi(x) := \exp(i(x,k)_H)$, $x \in H$, $h, k \in D(A^2)$. Then

$$\mathbb{E}[\varphi(Z_t(x))] = \exp(i(e^{-tA^2/2}h, x)_H - \tfrac{1}{2}((-A)^{-2}(I - e^{-tA^2})h, h)_H)$$

and computing the time derivative at $t = 0$ we obtain the generator of $Z$:

(2.8) $\qquad L\varphi(x) = -\frac{1}{2}\varphi(x)[i(A^2 h, x)_H + \|\Pi h\|_H^2].$

Now we compute the scalar product in $L^2(\mu_c; \mathbb{C})$ between $L\varphi$ and $\psi$:

$$\int L\varphi \, \bar{\psi} \, d\mu_c = -\frac{1}{2} \int [i(A^2 h, x)_H + \|\Pi h\|_H^2] \exp(i(h-k, x)_H) \mu_c(dx)$$

$$= -\frac{1}{2} \int [-(\Pi h, h-k)_H + \|\Pi h\|_H^2] \exp(i(h-k, x)_H) \mu_c(dx)$$

$$= -\frac{1}{2} \int (\Pi \nabla_H \varphi, \nabla_H \bar{\psi})_H \, d\mu_c,$$

where $\bar{\psi}$ is the complex conjugate of $\psi$ and in the second equality we have used (2.7). It follows that $(L, \operatorname{Exp}_A(H))$ is symmetric in $L^2(\mu_c)$ and the rest of the proof is standard. $\square$

**3. The approximating equation.** We begin with a few classical considerations on the approximating equation (1.5). First, notice that there is a conserved quantity, namely $t \mapsto \langle u_t^\varepsilon, e_0 \rangle_L$ is constant. Indeed, if we multiply the equation by $e_0 \equiv 1$ and we integrate in $\theta$, we obtain

$$\langle u_t^\varepsilon, e_0 \rangle_L = \langle u_0^\varepsilon, e_0 \rangle_L, \qquad t \geq 0.$$

In particular, if we want to study the ergodic properties of $u^\varepsilon$, we must restrict our attention to the affine subspace of initial conditions with fixed average $\{x \in L_c^2\}$.

Now we want to describe an important property of (1.5), namely, that it is a gradient system, see [10], Chapter 12. To this aim, it is convenient to write equation (1.5) in the abstract form:

(3.1) $\qquad \begin{cases} dX_t^\varepsilon = -\frac{1}{2} A(AX_t^\varepsilon - \nabla U_\varepsilon(X_t^\varepsilon)) \, dt + B \, dW_t, \\ X_0^\varepsilon(x) = x \in L^2(0, 1), \end{cases}$

with the notation already used in (2.2). Recall that $\nabla$ denotes the gradient in the Hilbert space $L^2(0, 1)$. Finally

$$U_\varepsilon(x) := \frac{1}{\varepsilon} \int_0^1 F(x(\theta)) \, d\theta, \qquad x \in L^2(0, 1), \qquad F'(u) = -f(u),$$

where we assume that $f : \mathbb{R} \mapsto \mathbb{R}$ is monotone nonincreasing and Lipschitz-continuous with $f(u) = 0$ for $u \geq 0$ and $f(u) > 0$ for $u < 0$. Notice that $\nabla U_\varepsilon(x) = -\frac{1}{\varepsilon} f(x)$ and $U_\varepsilon$ is a convex potential, since $-f$ is nondecreasing.

We write (3.1) in the mild formulation in $C([0, \infty); L^2(0, 1))$,

(3.2) $\qquad X_t^\varepsilon(x) = Z_t(x) - \frac{1}{2\varepsilon} \int_0^t A e^{-(t-s)A^2/2} f(X_s^\varepsilon) \, ds, \qquad t \geq 0,$

which is well defined since $\|A e^{-tA^2/2} h\|_L \leq \|h\|_L / \sqrt{t}$ for all $h \in L^2(0, 1)$. Then we have:

LEMMA 3.1. *Let $\varepsilon > 0$. For all $x \in H$ there exists a unique adapted process $X^\varepsilon \in C([0,\infty); L^2(0,1))$ solution of (3.1). Moreover, for all $t \geq 0$,*

$$\langle X_t^\varepsilon(x), e_0 \rangle_L = \langle x, e_0 \rangle_L. \tag{3.3}$$

PROOF. The proof of the first assertion is standard and based on a fixed point theorem. Indeed, since $f(\cdot)$ is Lipschitz-continuous, then the operator $\Gamma : C([0,T]; L^2(0,1)) \mapsto C([0,T]; L^2(0,1))$:

$$\Gamma(X)(t) := Z_t(x) - \frac{1}{2\varepsilon} \int_0^t A e^{-(t-s)A^2/2} f(X_s) \, ds, \qquad t \in [0,T],$$

is a contraction for the norm $\|X\|_\kappa := \sup_{t \in [0,T]} e^{-\kappa t} \|X_t\|_L$ for $\kappa > 0$ large enough. Formula (3.3) follows by taking the scalar product of both sides of (3.2) with $e_0$ in $L^2(0,1)$:

$$\langle X_t^\varepsilon(x), e_0 \rangle_L = \langle Z_t(x), e_0 \rangle_L = \langle x, e_0 \rangle_L$$

since $Ae_0 = 0$ and by (2.3). □

A crucial property of $X^\varepsilon$ is its 1-Lipschitz continuous dependence on the initial datum $x \in L_c^2$ in the norm of $H$ for every $c \in \mathbb{R}$: this is typical for dissipative systems. Notice however that this Lipschitz continuity fails if we want to let $x$ vary in $L^2(0,1)$ without the constraint of fixed average.

LEMMA 3.2. *Let $\varepsilon > 0$ and $c \in \mathbb{R}$. Then for all $t \geq 0$*

$$\|X_t^\varepsilon(x) - X_t^\varepsilon(y)\|_H \leq \exp(-t\pi^4/2) \|x - y\|_H, \qquad x, y \in L_c^2. \tag{3.4}$$

COROLLARY 3.3. *We can define by density a $H$-valued continuous process $(X_t^\varepsilon(x) : x \in H)$ which satisfies (3.4) for all $x, y \in H_c$ and solves (3.2) if $x \in L^2(0,1)$.*

REMARK 3.4. It is not difficult to prove that $X_t^\varepsilon(x)$ has paths in $L^4(0,T; L^2(0,1))$ and that it is in fact a solution of (3.1).

PROOF OF LEMMA 3.2. We consider for $N \in \mathbb{N}$ the process:

$$S_t^N := \sum_{i=0}^N \langle X_t^\varepsilon(x) - X_t^\varepsilon(y), e_i \rangle_L e_i, \qquad t \geq 0.$$

Then $t \mapsto S_t^N$ is $C^1$ with values in a $(N+1)$-dimensional subspace of $D(A)$. By (1.10) we have $-A\overline{Q}h = \Pi h$ for all $h \in L^2(0,1)$. Since $\Pi S_t^N = S_t^N$, by the spectral behavior of $A$ given in (1.8),

$$\frac{d}{dt} \|S_t^N\|_H^2 = \langle A S_t^N, S_t^N \rangle_L + \frac{1}{\varepsilon} \langle f(X_t^\varepsilon(x)) - f(X_t^\varepsilon(y)), S_t^N \rangle_L$$

$$\leq -\pi^4 \|S_t^N\|_H^2 + \frac{1}{\varepsilon} \langle f(X_t^\varepsilon(x)) - f(X_t^\varepsilon(y)), S_t^N \rangle_L.$$

This differential inequality implies:

$$\|S_t^N\|_H^2 \leq e^{-t\pi^4}\|x-y\|_H^2 + \int_0^t e^{-(t-s)\pi^4}\frac{1}{\varepsilon}\langle f(X_s^\varepsilon(x)) - f(X_s^\varepsilon(y)), S_s^N\rangle_L\,ds,$$

and by letting $N \to \infty$, since $f(\cdot)$ is monotone nonincreasing we obtain (3.4). □

We define for all $\varphi \in C_b(H_c)$ the semigroup and the resolvent of $X^\varepsilon$ on $H_c$:

$$P_t^{\varepsilon,c}\varphi(x) := \mathbb{E}[\varphi(X_t^\varepsilon(x))], \qquad x \in H_c, t \geq 0,$$

$$R_\lambda^{\varepsilon,c}\varphi(x) := \int_0^\infty e^{-\lambda t}\mathbb{E}[\varphi(X_t^\varepsilon(x))]\,dt, \qquad x \in H_c, \lambda > 0.$$

From (3.4) we deduce that $P_t^{\varepsilon,c}$ and $R_\lambda^{\varepsilon,c}$ act on $C_b(H_c)$ and moreover for all $c \in \mathbb{R}$ and $\varphi \in \mathrm{Lip}(H_c)$:

(3.5) $\quad \lambda|R_\lambda^{\varepsilon,c}\varphi(x) - R_\lambda^{\varepsilon,c}\varphi(y)| \leq [\varphi]_{\mathrm{Lip}}\|x-y\|_H, \qquad x,y \in H_c, \lambda > 0.$

We also define the probability measure on $L^2(0,1)$:

$$\nu_c^\varepsilon(dx) := \frac{1}{Z_c^\varepsilon}\exp(-U_\varepsilon(x))\mu_c(dx),$$

where $Z_c^\varepsilon$ is a normalization constant. Now, recalling (1.12) and (1.14), we introduce the symmetric bilinear form:

$$\mathcal{E}^{\varepsilon,c}(\varphi,\psi) := \tfrac{1}{2}\int_H (\Pi\nabla_H\varphi, \nabla_H\psi)_H\,d\nu_c^\varepsilon$$

$$= \tfrac{1}{2}\int_H \langle -A\nabla\varphi, \nabla\psi\rangle_L\,d\nu_c^\varepsilon \qquad \forall\,\varphi,\psi \in \mathrm{Exp}_A(H).$$

In the following result we prove that $X^\varepsilon$ is strong Feller in $H_c$ for all $c \in \mathbb{R}$ and is associated with (the closure of) $\mathcal{E}^{\varepsilon,c}$, that $\{\nu_c^\varepsilon : c \in \mathbb{R}\}$ is the set of all ergodic invariant probability measures of $X^\varepsilon$. Moreover we prove that $X^\varepsilon$ is reversible with respect to each $\nu_c^\varepsilon$ for $c \in \mathbb{R}$.

PROPOSITION 3.5. *Let $c \in \mathbb{R}$ and $\varepsilon > 0$.*

1. $(\mathcal{E}^{\varepsilon,c}, \mathrm{Exp}_A(H))$ *is closable in* $L^2(\nu_c^\varepsilon)$: *we denote by* $(\mathcal{E}^{\varepsilon,c}, D(\mathcal{E}^{\varepsilon,c}))$ *the closure. Moreover* $\mathrm{Lip}(H_c) \subset D(\mathcal{E}^{\varepsilon,c})$ *and* $\mathcal{E}^{\varepsilon,c}(\varphi,\varphi) \leq [\varphi]_{\mathrm{Lip}(H_c)}^2$.
2. $(R_\lambda^{\varepsilon,c})_{\lambda>0}$ *is the resolvent associated with* $\mathcal{E}^{\varepsilon,c}$, *that is, for all* $\lambda > 0$ *and* $\varphi \in L^2(\nu_c^\varepsilon)$, $R_\lambda^{\varepsilon,c}\varphi \in D(\mathcal{E}^{\varepsilon,c})$ *and*

(3.6) $\quad \lambda\int_H R_\lambda^{\varepsilon,c}\varphi\psi\,d\nu_c^\varepsilon + \mathcal{E}^{\varepsilon,c}(R_\lambda^{\varepsilon,c}\varphi, \psi) = \int_H \varphi\psi\,d\nu_c^\varepsilon \qquad \forall\,\psi \in D(\mathcal{E}^{\varepsilon,c}).$

3. $\nu_c^\varepsilon$ is an ergodic invariant probability measure of $X^\varepsilon$ and $X^\varepsilon$ is reversible with respect to $\nu_c^\varepsilon$. Moreover for all $\varphi \in \mathrm{Lip}(H_c)$:

$$(3.7) \qquad \lim_{t\to\infty} |\mathbb{E}[\varphi(X_t^\varepsilon(x))] - \nu_c^\varepsilon(\varphi)| = 0, \qquad x \in H_c,$$

and $\{\nu_c^\varepsilon : c \in \mathbb{R}\}$ are the only ergodic invariant probability measures of $X^\varepsilon$.

4. For all $\varphi : H_c \mapsto \mathbb{R}$ bounded and Borel we have

$$(3.8) \qquad |P_t^{\varepsilon,c}\varphi(x) - P_t^{\varepsilon,c}\varphi(y)| \leq \frac{\|\varphi\|_\infty}{\sqrt{t}} \|x - y\|_H, \qquad x, y \in H_c,\ t > 0.$$

In particular, $X^\varepsilon$ is strong Feller on $H_c$.

PROOF. The proof of points 1 and 2 is standard, see [20] and [10], Chapter 12, so we only sketch the proof. By (3.2) the process $t \mapsto (X_t^\varepsilon(x), h)_H$ is a semimartingale for $h \in D(A)$ and for $t \geq 0$

$$(X_t^\varepsilon(x), h)_H = (x, h)_H + \frac{1}{2} \int_0^t \langle X_s^\varepsilon(x), Ah \rangle_L\, ds$$
$$+ \frac{1}{2\varepsilon} \int_0^t \langle f(X_s^\varepsilon(x)), \Pi h \rangle_L\, ds + M_t^h,$$

where $M^h$ is a martingale with quadratic variation $\langle M^h \rangle_t = t \|\Pi h\|_H^2$. Like in the proof of Proposition 2.2, let $\varphi(x) := \exp(i(x, h)_H)$ for $x \in H$ and $h \in D(A^2)$. By Itô's formula,

$$(3.9) \quad L^\varepsilon \varphi(x) := \frac{d}{dt} \mathbb{E}[\varphi(X_t^\varepsilon(x))]\bigg|_{t=0} = L\varphi(x) + \frac{i}{2\varepsilon} \langle f(x), \Pi h \rangle_L \varphi(x),$$

where $L\varphi$ is as in (2.8). An application of (2.7) shows that $(L^\varepsilon, \mathrm{Exp}_A(H))$ is symmetric in $L^2(\nu_c^\varepsilon)$ and

$$\int L^\varepsilon \varphi \psi\, d\nu_c^\varepsilon = -\tfrac{1}{2} \int (\Pi \nabla_H \varphi, \nabla_H \psi)_H\, d\nu_c^\varepsilon \qquad \forall\, \varphi, \psi \in \mathrm{Exp}_A(H),$$

and the rest of the proof is standard.

From (3.6) we obtain in particular symmetry of $R_\lambda^{\varepsilon,c}$ in $L^2(\nu_c^\varepsilon)$, hence reversibility of $X^\varepsilon$ w.r.t. $\nu_c^\varepsilon$. Now, the proof of (3.7) is based on a standard coupling method and on Lemma 3.2: let $\mathcal{Y}$ be a $H$-valued r.v. with distribution $\nu_c^\varepsilon$ and independent of $W$. Then by (3.4)

$$|\mathbb{E}[\varphi(X_t^\varepsilon(x))] - \nu_c^\varepsilon(\varphi)| = |\mathbb{E}[\varphi(X_t^\varepsilon(x)) - \varphi(X_t^\varepsilon(\mathcal{Y}))]|$$
$$\leq 2\mathbb{E}[\|\varphi\|_\infty \wedge ([\varphi]_{\mathrm{Lip}(H_c)} \|x - \mathcal{Y}\|_H \exp(-t\pi^4/2))]$$

and the thesis follows by dominated convergence.

Now we prove the strong Feller property. Fix $c \in \mathbb{R}$. We notice that the process $H_0 \ni x \mapsto \mathcal{X}_t(x) := X_t^\varepsilon(x + ce_0) - ce_0 \in H_0$ solves the following equation:

$$(3.10) \qquad \begin{cases} d\mathcal{X}_t = -\frac{1}{2} A(A\mathcal{X}_t - \nabla U_\varepsilon(ce_0 + \mathcal{X}_t))\, dt + B\, dW_t, \\ \mathcal{X}_0(x) = x \in H_0. \end{cases}$$

This process is a gradient system in $H_0$ with nondegenerate noise and with a convex potential $U_\varepsilon(ce_0 + \cdot)$: the proof of the strong Feller property can be found in [10], Chapter 12.3; see also [6]. □

**4. Pathwise uniqueness of solutions of (1.1).** In this section we turn our attention to equation (1.1) and we prove that for any pair $(u^i, \eta^i, W)$, $i = 1, 2$, of weak solutions of (1.1) defined on the same probability space with the same driving noise $W$ and with $u_0^1 = u_0^2$, we have $(u^1, \eta^1) = (u^2, \eta^2)$. This pathwise uniqueness is used in the next section to construct stationary strong solutions of (1.1).

PROPOSITION 4.1. *Let $u_0 \in C([0,1])$ with $u_0 \geq 0$, $\int_0^1 u > 0$ and $(u_0, W)$ independent. Let $(u^1, \eta^1, W)$ and $(u^2, \eta^2, W)$ be two weak solutions of (1.1) with $u_0^1 = u_0^2$. Then $(u^1, \eta^1) = (u^2, \eta^2)$.*

Throughout the paper we use several times the following easy result:

LEMMA 4.2. *Let $\zeta(dt, d\theta)$ be a finite signed measure on $[\delta, T] \times [0, 1]$ and $v \in C([\delta, T] \times [0, 1])$. Suppose that for all $s \in [\delta, T]$*

$$(4.1) \qquad \int_{[s,T] \times [0,1]} h_\theta \zeta(dt, d\theta) = 0 \qquad \forall\, h \in C([0,1]), \bar{h} = 0,$$

*and*

$$(4.2) \qquad \bar{v}_s = c > 0, \qquad \int_{[s,T] \times [0,1]} v\, d\zeta = 0.$$

*Then $\zeta \equiv 0$.*

PROOF. Setting $h := k - \bar{k}$, $k \in C([0,1])$, we obtain by (4.1) for all $\delta \leq s \leq t \leq T$

$$\int_0^1 k_\theta \zeta([s,t] \times d\theta) = \zeta([s,t] \times [0,1]) \int_0^1 k_\theta\, d\theta \qquad \forall\, k \in C([0,1]).$$

This implies $\zeta(dt, d\theta) = \gamma(dt)\, d\theta$, where $\gamma(t) := \zeta([\delta, t] \times [0, 1])$, $t \in [\delta, T]$, is a process with bounded variation. Then by (4.2)

$$0 = \int_{[s,t] \times [0,1]} v\, d\zeta = \int_s^t \left( \int_0^1 v_s(\theta)\, d\theta \right) \gamma(ds) = c(\gamma(t) - \gamma(s)),$$

that is, $\gamma(t) - \gamma(s) = 0$, since $c > 0$. □

PROOF OF PROPOSITION 4.1. Let us set $v := u^1 - u^2$, $\zeta := \eta^1 - \eta^2$. By (1.13), for all $h \in D(A^2)$ and $t \geq \delta > 0$,

$$(4.3) \quad \langle v_t, h \rangle_L = \langle v_\delta, h \rangle_L - \tfrac{1}{2} \int_\delta^t \langle v_s, A^2 h \rangle_L\, ds - \tfrac{1}{2} \int_{[\delta, t] \times [0,1]} Ah_\theta \zeta(ds, d\theta).$$

We consider the following approximation of $v$:
$$v_t^N := \frac{1}{N} \sum_{n=0}^{N} \sum_{i=0}^{n} \langle v_t, e_i \rangle_L e_i.$$

Since $v$ is continuous, then $v^N$ converges uniformly to $v$ on $[0,T] \times [0,1]$. Notice that for all $i \geq 0$, the process $[\delta, T] \ni t \mapsto \langle v_t, e_i \rangle_L$ has bounded variation by (4.3). In particular, the process $[\delta, T] \ni t \mapsto v_t^N$ has bounded variation with values in a finite-dimensional subspace of $D(A)$. By Itô's formula,

$$\begin{aligned}(v_t, v_t^N)_H &= \frac{1}{N} \sum_{n=0}^{N} \sum_{i=0}^{n} (i\pi)^{-2} \langle v_t, e_i \rangle_L^2 \\ &= (v_\delta, v_\delta^N)_H - \int_\delta^t \frac{1}{N} \sum_{n=0}^{N} \sum_{i=0}^{n} (i\pi)^2 \langle v_s, e_i \rangle_L^2 \, ds \\ &\quad + \int_{[\delta,t] \times [0,1]} v_s^N(\theta) \zeta(ds, d\theta) \\ &\leq (v_\delta, v_\delta^N)_H + \int_{[\delta,t] \times [0,1]} v_s^N(\theta) \zeta(ds, d\theta).\end{aligned}$$

Letting $N \to \infty$, $v^N$ converges uniformly on $[0,t] \times [0,1]$ to $v$, which is continuous. By dominated convergence,

$$\|v_t\|_H^2 - \|v_\delta\|_H^2 \leq \int_{[\delta,t] \times [0,1]} v \, d\zeta = -\int_{[\delta,t] \times [0,1]} u^1 \, d\eta^2 - \int_{[\delta,t] \times [0,1]} u^2 \, d\eta^1 \leq 0,$$

and letting $\delta \to 0$ we obtain $v = 0$ and $u^1 = u^2$. Turning to $\zeta$, we see by (4.3) that $\int_{[\delta,t] \times [0,1]} Ah_\theta \zeta(dt, d\theta) = 0$ for all $h \in D(A^2)$. By density, we obtain that $\zeta$ and $u_1 = u_2 =: v$ satisfy (4.1) and (4.2) above, and therefore by Lemma 4.2, $\zeta \equiv 0$, that is, $\eta^1 = \eta^2$. □

**5. Existence of stationary solutions of (1.1).** In this section we prove the existence of stationary strong solutions of equation (1.1) and that they are limits in distribution of stationary solutions of (1.5). Fix $c > 0$ and consider the unique (in law) stationary solution of (1.5), $(\widehat{X}^{\varepsilon,c})$, in $H_c$. We are going to prove that the laws of $(\widehat{X}^{\varepsilon,c})_{\varepsilon>0}$ weakly converge as $\varepsilon \to 0$ to a stationary weak solution of (1.1).

THEOREM 5.1. *Let $c > 0$. For all $T > 0$, $\widehat{X}^{\varepsilon,c}$ weakly converges in $C([0,T] \times [0,1])$ to a process $\widehat{X}^c$ as $\varepsilon \to 0$. Moreover setting $u^\varepsilon := \widehat{X}^{\varepsilon,c}$ and*
$$\eta^\varepsilon(dt, d\theta) := \frac{1}{\varepsilon} f(u_t^\varepsilon(\theta)) \, dt \, d\theta,$$
*then $(u^\varepsilon, \eta^\varepsilon, W)$ converges in law to $(u, \eta, W)$, stationary weak solution of (1.1).*

The proof is based on two steps: tightness and identification of the limit. First, we set:

(5.1) $$K := \{x \in L^2(0,1) : x \geq 0\}.$$

By Lemma 2.1, $\mu_c$ is the distribution of $Y^c := B - \overline{B} + c$. Now, notice the following inclusion of events:

$$\{B_\theta \in [-c/2, c/2], \forall \theta \in [0,1]\} \subset \{Y_c \in K\}.$$

Therefore $\mu_c(K) > 0$ for $c > 0$ and, as $\varepsilon \to 0$,

(5.2) $\quad c > 0 \quad \Longrightarrow \quad \nu_c^\varepsilon \rightharpoonup \nu_c(dx) := \mu_c(dx|K) = \dfrac{1}{\mu_c(K)} 1_{(x \in K)} \mu_c(dx).$

In particular, the initial distribution of $\widehat{X}^{\varepsilon,c}$ converges to $\nu_c$.

LEMMA 5.2. *For all $T > 0$, the laws of $(\widehat{X}^{\varepsilon,c})_{\varepsilon > 0}$ are tight in $C([0,T] \times [0,1])$.*

PROOF. We introduce the space $H^{-\gamma}(0,1)$, $\gamma > 0$, completion of $L^2(0,1)$ w.r.t. the norm

$$\|f\|_{-\gamma}^2 := \sum_{n=0}^\infty (1+n)^{-2\gamma} |\langle f, e_n \rangle_L|^2$$

where $e_n$ is defined in (1.8). Notice that $H^{-1}(0,1) = H$, in our notation. We recall that the Hilbert–Schmidt norm of the inclusion $H = H^{-1}(0,1) \to H^{-\gamma}(0,1)$ is finite for all $\gamma > 2$. We fix $\gamma > 2$. We claim that for all $p > 1$ there exists $C_p \in (0, \infty)$, independent of $\varepsilon$, such that for all $N \in \mathbb{N}$:

(5.3) $\quad (\mathbb{E}[\|\widehat{X}_t^{\varepsilon,c} - \widehat{X}_s^{\varepsilon,c}\|_{H^{-\gamma}(0,1)}^p])^{1/p} \leq C_p |t-s|^{1/2}, \qquad t, s \geq 0.$

To prove (5.3), we fix $\varepsilon > 0$ and $T > 0$ and use the Lyons–Zheng decomposition (see, e.g., [15], Theorem 5.7.1) to write for $t \in [0, T]$ and $h \in H$:

$$(h, \widehat{X}_t^{\varepsilon,c} - \widehat{X}_0^{\varepsilon,c})_H = \tfrac{1}{2} M_t - \tfrac{1}{2}(N_T - N_{T-t}),$$

where $M$, respectively, $N$, is a martingale w.r.t. the natural filtration of $\widehat{X}^{\varepsilon,c}$, respectively, of $(\widehat{X}_{T-t}^{\varepsilon,c}, t \in [0,T])$. Moreover, the quadratic variations are both equal to: $\langle M \rangle_t = \langle N \rangle_t = t \cdot \|\Pi h\|_H^2$. By the Burkholder–Davis–Gundy inequality we can find $c_p \in (0, \infty)$ for all $p > 1$ such that:

$$(\mathbb{E}[\|\widehat{X}_t^{\varepsilon,c} - \widehat{X}_s^{\varepsilon,c}\|_{H^{-\gamma}(0,1)}^p])^{1/p} \leq c_p \kappa_{-\gamma} |t-s|^{1/2}, \qquad t, s \in [0, T],$$

where $\kappa_{-\gamma}$ is the Hilbert–Schmidt norm of the inclusion of $H^{-1}(0,1)$ into $H^{-\gamma}(0,1)$.

We also have for $\eta < 1/2$ and $r \geq 1$ by stationarity

$$(\mathbb{E}[\|\widehat{X}_t^{\varepsilon,c} - \widehat{X}_s^{\varepsilon,c}\|_{W^{\eta,r}(0,1)}^p])^{1/p}$$

(5.4)
$$\leq (\mathbb{E}[\|\widehat{X}_t^{\varepsilon,c}\|_{W^{\eta,r}(0,1)}^p])^{1/p} + (\mathbb{E}[\|\widehat{X}_s^{\varepsilon,c}\|_{W^{\eta,r}(0,1)}^p])^{1/p}$$

$$= 2\left(\int_H \|x\|_{W^{\eta,r}(0,1)}^p d\nu_c^\varepsilon\right)^{1/p} \leq c\left(\int_H \|x\|_{W^{\eta,r}(0,1)}^p d\mu_c\right)^{1/p}$$

since $U_\varepsilon \geq 0$. The latter term is finite by the representation of Lemma 2.1. Let us now take $\kappa \in [0,1]$ and set $\alpha = \kappa\eta - (1-\kappa)\gamma$, $\frac{1}{q} = \kappa\frac{1}{r} + (1-\kappa)\frac{1}{2}$, then by interpolation

$$(\mathbb{E}[\|\widehat{X}_t^{\varepsilon,c} - \widehat{X}_s^{\varepsilon,c}\|_{W^{\alpha,q}(0,1)}^p])^{1/p}$$

$$\leq (\mathbb{E}[\|\widehat{X}_t^{\varepsilon,c} - \widehat{X}_s^{\varepsilon,c}\|_{W^{\eta,r}(0,1)}^p])^{\kappa/p}(\mathbb{E}[\|\widehat{X}_t^{\varepsilon,c} - \widehat{X}_s^{\varepsilon,c}\|_{H^{-\gamma}(0,1)}^p])^{(1-\kappa)/p}.$$

For any $\beta \in (0, 1/2)$, we can choose $\eta \in (0, 1/2)$, $\gamma > 2$, $r \geq 1$ and $\kappa \in (0,1)$ such that $(\alpha - \beta)q > 1$. It follows, by the Sobolev embedding and (5.3) and (5.4), that

$$(\mathbb{E}[\|\widehat{X}_t^{\varepsilon,c} - \widehat{X}_s^{\varepsilon,c}\|_{C^\beta([0,1])}^p])^{1/p} \leq \tilde{c}|t-s|^{(1-\kappa)/(2p)}.$$

Since the law of $\widehat{X}_0^{\varepsilon,c}$ is $\nu_c^\varepsilon$ which converges as $\varepsilon \to 0$, tightness of the laws of $(\widehat{X}^{\varepsilon,c})_{\varepsilon>0}$ in $C([0,T] \times [0,1])$ follows, for example, by Theorem 7.2 in Chapter 3 of [14]. $\square$

We define the Polish space $E := C(O_T) \times M(O_T) \times C(O_T)$, where $O_T := [0,T] \times [0,1]$ and $M(O_T)$ is the space of all finite positive measures on $[0,T] \times [0,1]$ endowed with the weak topology of the dual space of $C(O_T)$.

LEMMA 5.3. *Let $c > 0$. Let $\varepsilon_n \downarrow 0$ be any sequence such that $u^{\varepsilon_n}$ converges in law to a process $u$. Then $(u^{\varepsilon_n}, \eta^{\varepsilon_n}, W)$ converges in law to $(u, \eta, W)$, stationary weak solution of (1.1), in $E$.*

PROOF. By Skorohod's theorem we can find a probability space and a sequence of processes $(v_n, w_n)$ such that $(v_n, w_n) \to (v, w)$ in $C(O_T)$ almost surely and $(v_n, w_n)$ has the same distribution as $(u^{\varepsilon_n}, W)$ for all $n \in \mathbb{N}$, where $O_T := [0,T] \times [0,1]$. Notice that $v \geq 0$ almost surely, since for all $t$ the law of $v_t(\cdot)$ is $\nu_c$ which is concentrated on $K$ and moreover $v$ is continuous on $O_T$. We set now

$$\xi^n(dt, d\theta) := \frac{1}{\varepsilon_n} f(v_t^n(\theta)) dt d\theta.$$

From (1.5) we obtain that a.s. for all $T \geq 0$ and $h \in D(A^2)$ and $\bar{h} = 0$,

$$(5.5) \qquad \exists \lim_{n \to \infty} \int_{O_T} h_\theta \xi^n(dt, d\theta).$$

The limit is a random distribution on $O_T$. We want to prove that in fact $\xi^n$ converges as a measure in the dual of $C(O_T)$ for all $T \geq 0$. For this, it is enough to prove that the mass $\xi^n(O_T)$ converges as $n \to \infty$.

Suppose that $\{\xi^n(O_T)\}_n$ is unbounded. We define $\zeta^n := \xi^n/\xi^n(O_T)$. Then $\zeta^n$ is a probability measure on the compact set $O_T$. By tightness we can extract from any subsequence a sub-subsequence converging to a probability measure $\zeta$. By the uniform convergence of $v^n$ we can see that the contact condition $\int_{O_T} v \, d\zeta = 0$ holds. Moreover, dividing (1.5) by $\xi^n(O_T)$ for $t \in [0,T]$, we obtain that $\int_{O_t} h_\theta \zeta(ds, d\theta) = 0$ for all $h \in D(A^2)$ with $\bar{h} = 0$ and by density for all $h \in C([0,1])$ with $\bar{h} = 0$.

Then $\zeta$ and $v$ satisfy (4.1) and (4.2) above, and therefore by Lemma 4.2, $\zeta \equiv 0$, a contradiction since $\zeta$ is a probability measure. Therefore $\limsup_{n \to \infty} \xi^n(O_T) < \infty$.

By tightness, for any subsequence in $\mathbb{N}$ we have convergence of $\xi^n$ to a finite measure $\xi$ on $[0,T] \times [0,1]$ along some sub-subsequence. Let $\xi_i$, $i = 1, 2$, be two such limits and set $\zeta := \xi_1 - \xi_2$. By (5.5) and by density

$$\int_{O_T} h_\theta \xi_1(dt, d\theta) = \int_{O_T} h_\theta \xi_2(dt, d\theta) \qquad \forall \, h \in C([0,1]), \bar{h} = 0,$$

that is $\zeta$ and $v$ satisfy (4.1) and (4.2) above. By Lemma 4.2, $\zeta \equiv 0$, that is, $\xi_1 = \xi_2$. Therefore, $\xi_n$ converges as $n \to \infty$ to a finite measure $\xi$ on $[0,T] \times [0,1]$.

Finally, we need to prove that the contact condition holds, that is, that $\int_{(0,\infty) \times [0,1]} v \, d\xi = 0$. Since $f \geq 0$ and $f(u) > 0$ for $u > 0$, then $uf(u) \leq 0$ for all $u \in \mathbb{R}$. Then

$$0 \geq \int_{[0,T] \times [0,1]} v^n \, d\xi^n \to \int_{[0,T] \times [0,1]} v \, d\xi$$

by the uniform convergence of $v^n$ to $v$, the convergence of $\xi^n$ to $\xi$ and Lemma 8.2 below. Since $v \geq 0$ and $\xi$ is a positive measure, then $\int_{[0,T] \times [0,1]} v \, d\xi \leq 0$ is possible only if $\int_{[0,T] \times [0,1]} v \, d\xi = 0$ $\square$

LEMMA 5.4. *Let $c > 0$. As $\varepsilon \to 0$, $u^\varepsilon$ converges in law. Moreover, every weak stationary solution of (1.1) with $\bar{u}_0 = c$ a.s. is also a strong solution.*

PROOF. We use a technique presented in [19]. It is no loss of generality to assume that $u_0^\varepsilon(\cdot)$ converges in probability as $\varepsilon \to 0$. Let $(\varepsilon_n^1)_{n \in \mathbb{N}}$ and $(\varepsilon_n^2)_{n \in \mathbb{N}}$ be two sequences of positive numbers converging to 0 as $n \to \infty$. In the notation of Lemma 5.3, by Lemma 5.2 the process $(u^{\varepsilon_n^1}, u^{\varepsilon_n^2}, W)$ is tight in a suitable space. By Skorohod's theorem we can find a probability space and a sequence of processes $(v_n^1, v_n^2, w_n)$ such that $(v_n^1, v_n^2, w_n) \to (v^1, v^2, w)$ a.s. and $(v_n^1, v_n^2, w_n)$ has the same distribution as $(u^{\varepsilon_n^1}, u^{\varepsilon_n^2}, W)$ for all $n \in \mathbb{N}$.

By Lemma 5.3, $(v_1, w)$ and $(v_2, w)$ are both weak solutions. By Proposition 4.1, necessarily $v_1 = v_2$. Therefore the process $u^{\varepsilon_n^1} - u^{\varepsilon_n^2}$ converges in law to the process constantly zero, and therefore it converges in probability. It follows that the sequence $(u^\varepsilon)$ is Cauchy and converges in probability and therefore in law to $u$ which is a stationary strong solution of (1.1).

By pathwise uniqueness and existence of strong solutions, we obtain that every weak solution is also strong. □

The last three lemmas yield the proof of Theorem 5.1. Moreover, we have the following result:

COROLLARY 5.5. *Let $c > 0$.*

1. *There exists a continuous process $(X_t(x), t \geq 0, x \in K \cap H_c)$ with $X_0(x) = x$ and a set $K_0$ dense in $K \cap H_c$, such that for all $x \in K_0$ there exists a unique strong solution $(u, \eta, W)$ of (1.1) with $u_t = X_t(x)$, $t \geq 0$.*
2. *The law of $(X_t(x), t \geq 0)$ is a regular conditional distribution of the law of $\widehat{X}^c$ given $\widehat{X}_0^c = x \in K \cap H_c$.*

PROOF. By Lemma 5.4, we have a stationary strong solution $u$ in $H_c$ with $W$ and $u_0$ independent. Conditioning on the value of $u_0$, we obtain for $\nu_c$-a.e. $x$ a strong solution $u$ with $u_0 = x$. Since the support of $\nu_c$ is $K \cap H_c$, we have a strong solution for a dense set $K_0$ of $x$ in $K \cap H_c$.

Notice that all processes $(X_t(x), t \geq 0)$ with $x \in K_0$ are driven by the same noise $W$ and are in particular continuous with values in $H$. Arguing like in the proof of Proposition 4.1, we see that

$$(5.6) \qquad \|X_t(x) - X_t(y)\|_H \leq \|x - y\|_H \qquad \forall\, x, y \in K_0,\ t \geq 0.$$

Then, by density, we obtain a continuous process $(X_t(x) : t \geq 0)$ in $H_c$ for all $x \in K \cap H_c$. □

Notice that, in Corollary 5.5, we are not able yet to say that $(X_t(x), t \geq 0)$ is a solution (and therefore the unique one) of (1.1) for $x \notin K_0$. Indeed, the equation requires that the solution be continuous on $(0, T] \times [0, 1]$ a.s. for the contact condition (1.2) to be meaningful.

This problem is solved in the next section. The crucial point will be the strong Feller property for the transition semigroup of $(X_t, t \geq 0)$; see Proposition 6.1.

**6. Existence of solutions of (1.1).** We want now to prove that for any deterministic initial condition $u_0 = x \in K \cap H_c$, $c > 0$, there exists a strong solution of equation (1.1), necessarily unique by Proposition 4.1, and that the process $X$ constructed in Corollary 5.5 is a realization of such solution.

We recall that in Corollary 5.5 we constructed a continuous process $(X_t(x), t \geq 0, x \in K \cap H_c)$ such that $(X_t(x), t \geq 0)$ is a strong solution of (1.1) for $x$ in a dense set. We prove now that the transition semigroup $P^{\varepsilon,c}$ of $X^\varepsilon$ on $K \cap H_c$ converges to the transition semigroup $P^c$ of $X$ on $K \cap H_c$. The result is the following:

PROPOSITION 6.1. *Let $c > 0$. For all $\varphi \in C_b(H)$ and $x \in K \cap H_c$:*

$$(6.1) \qquad \lim_{\varepsilon \to 0} P_t^{\varepsilon,c} \varphi(x) = \mathbb{E}[\varphi(X_t(x))] =: P_t^c \varphi(x).$$

*Moreover the Markov process $(X_t(x), t \geq 0, x \in K \cap H_c)$ is strong Feller:*

$$(6.2) \quad |P_t^c \varphi(x) - P_t^c \varphi(y)| \leq \frac{\|\varphi\|_\infty}{\sqrt{t}} \|x - y\|_H, \qquad x, y \in K \cap H_c,\ t > 0.$$

PROOF. Fix $t > 0$. By (3.8), for any $\varphi \in C_b(H) : \sup_\varepsilon (\|P_t^{\varepsilon,c} \varphi\|_\infty + [P_t^{\varepsilon,c} \varphi]_{\mathrm{Lip}(H_c)}) < \infty$. Let $(\varepsilon_j)_j$ be any sequence in $\mathbb{N}$ and $(x_k)_k$ a countable dense set in $H_c$. With a diagonal procedure, by Ascoli–Arzelà's theorem we can find a subsequence $(j_i)_i$ and a function $F : \{x_k, k \in \mathbb{N}\} \mapsto \mathbb{R}$ such that $P_t^{\varepsilon_{j_i},c} \varphi(x_k) \to F(x_k)$ as $j = j_i \to \infty$ for all $k \in \mathbb{N}$. By (3.5), $F$ is Lipschitz on $\{x_k, k \in \mathbb{N}\}$ and therefore can be extended to a function in $F \in \mathrm{Lip}(H_c)$ and

$$(6.3) \qquad F(x) = \lim_{i \to \infty} P_t^{\varepsilon_{j_i},c} \varphi(x) \qquad \forall\, x \in H_c.$$

On the other hand, by Theorem 5.1, for all $\varphi \in C_b(H)$,

$$\mathbb{E}[\psi(\widehat{X}_0^c)\varphi(\widehat{X}_t^c)] = \lim_{i \to \infty} \mathbb{E}[\psi(\widehat{X}_0^{\varepsilon_{j_i},c})\varphi(\widehat{X}_t^{\varepsilon_{j_i},c})]$$
$$= \lim_{i \to \infty} \int \psi P_t^{\varepsilon_{j_i},c} \varphi\, d\nu_c^{\varepsilon_{j_i}} = \int \psi F\, d\nu_c.$$

By Corollary 5.5 we obtain: $F(x) = \mathbb{E}[\varphi(X_t(x))]$ for $\nu_c$-a.e. $x$. By (5.6), both sides of the latter equality are continuous in $x \in K \cap H_c$ and therefore they coincide for all $x \in K \cap H_c$. Then the limit in (6.3) does not depend on the chosen subsequence $(\varepsilon_{j_i})_i$ and we obtain (6.1). By (3.8) we obtain (6.2). □

THEOREM 6.2. *Let $u_0$ be a $K$-valued r.v. with $\bar{u}_0 > 0$ a.s. and $(u_0, W)$ independent. Then:*

1. *There exists a unique strong solution $(u, \eta, W)$ of (1.1); moreover $X_t(u_0) = u_t$, $t \geq 0$.*

2. *The process $(X_t(x) : t \geq 0, x \in K \cap H_c)$ is a Markov process with transition semigroup $P^c$, continuous and strong Feller in $H_c$.*
3. *For all $c > 0$, $x \in K \cap H_c$ and $0 = t_0 < t_1 < \cdots < t_n$, $(X_{t_i}(x), i = 1, \ldots, n)$ is the limit in distribution of $(X^\varepsilon_{t_i}(x), i = 1, \ldots, n)$.*
4. *If $u_0$ has distribution $\nu_c$, $c > 0$, then $(X_t(u_0), t \geq 0)$ is equal in distribution to $(\widehat{X}^c_t, t \geq 0)$, see Theorem 5.1.*

In particular, $(X_t(x) : t \geq 0, x \in K \cap H_c)$ is the limit of $(X^\varepsilon_t(x) : t \geq 0, x \in K \cap H_c)$ in the sense of the finite-dimensional distributions.

PROOF OF THEOREM 6.2. In Theorem 5.1 we proved convergence of $\widehat{X}^{\varepsilon,c}$ to $\widehat{X}^c$. By Corollary 5.5 we have a process $(X_t(x), t \geq 0, x \in H \cap H_c)$, such that for all $x$ in a set $K_0$ dense in $H \cap H_c$ we have a strong solution of (1.1) with initial condition $u_0 = x$ and with $u_t = X_t(x)$ for all $t \geq 0$. By Proposition 6.1 we have now that the Markov process $X$ has transition semigroup $P^c$ on $H_c$.

Notice now that the strong Feller property (6.2) of $P^c_t$ implies that for all $x \in K \cap H_c$ and $s > 0$ the law of $X_s(x)$ is absolutely continuous w.r.t. the invariant measure $\nu_c$. Indeed, if $\nu_c(\Gamma) = 0$, then $\nu_c(P^c_s(1_\Gamma)) = \nu_c(\Gamma) = 0$ so that $P^c_s(1_\Gamma)(x) = 0$ for $\nu_c$-a.e. $x$ and by continuity for all $x \in K \cap H_c$.

Therefore, a.s. $X_s(x) \in K_0$ for all $s > 0$ and $x \in K \cap H_c$ and in particular $(X_{t+s}(x), t \geq 0)$ is a strong solution of (1.1) with initial condition $X_s(x)$. In particular, we have a process $u \in C([0,T]; H) \cap C((0,T] \times [0,1])$ and a measure $\eta$ on $(0,T] \times [0,1]$ which is finite on $[\delta, T] \times [0,1]$ for all $\delta > 0$, such that $(u_{t+s}, \eta(s + dt, d\theta), W(t+s, \cdot) - W(s, \cdot))_{t \geq 0}$ is a strong solution of (1.1) with initial condition $X_s(x)$ for all $s > 0$. Since $X_s(x) \to x$ in $H$ as $s \to 0$, then $(u, \eta, W)$ is a strong solution of (1.1) with initial condition $u_0 = x$ in the sense of Definition 1.1.

By Proposition 6.1, the process $(X_t(x) : t \geq 0, x \in K \cap H_c)$ has transition semigroup $P^c$ and is the limit of $(X^\varepsilon_t(x) : t \geq 0, x \in K \cap H_c)$ in the sense of the finite-dimensional distributions. $\square$

REMARK 6.3. We have not been able to prove that $\eta$ is finite on $[0, T] \times [0, 1]$ nor that $u$ is continuous on $[0, T] \times [0, 1]$. Both properties are true for stationary solutions and for $\nu_c$-a.e. initial condition $x$; see Lemma 5.3 and Corollary 5.5. However, notice that for all $h \in D(A^2)$, by (1.13) and by the continuous dependence of $x \mapsto X(x)$ in $H$, we can see that for all initial condition $x \in C([0, 1]) \cap K$ and $t > 0$,

$$\exists \lim_{\delta \downarrow 0} \int_\delta^t \int_0^1 Ah_\theta \eta(ds, d\theta) \in \mathbb{R}.$$

However, since $\overline{Ah} = 0$, we have no information on the mass $\eta([\delta, T] \times [0, 1])$ as $\delta \to 0$,

We also have the following result concerning the ergodic properties of $X$:

PROPOSITION 6.4. *For all $c > 0$, $\nu_c$ is an ergodic invariant probability measure of $X$ and $X$ is reversible with respect to $\nu_c$. Moreover for all $\varphi \in \mathrm{Lip}(H_c)$,*

$$\lim_{t \to \infty} |\mathbb{E}[\varphi(X_t(x))] - \nu_c(\varphi)| = 0, \qquad x \in H_c, \tag{6.4}$$

*and $\{\nu_c : c \in \mathbb{R}\}$ are the only ergodic invariant probability measures of $X$,*

This result can be proven like point 3 of Proposition 3.5.

**7. An integration by parts formula.** From now on we consider $c > 0$. We want to prove an integration by parts formula on the infinite-dimensional probability measure $\nu_c$, defined in (5.2) above. We recall that, by Lemma 2.1, if $B$ a Brownian motion and $c > 0$ is constant, then $\nu_c$ is the law of the process $Y_\theta^c := B_\theta - \overline{B} + c$ conditioned to be nonnegative on $[0, 1]$.

We denote by $(M, \hat{M})$ two independent copies of the standard Brownian meander (see [22]) and we set, for all $r \in (0, 1)$,

$$U_r(\theta) := \begin{cases} \sqrt{r} M\left(\dfrac{r - \theta}{r}\right), & \theta \in [0, r], \\ \sqrt{1 - r} \hat{M}\left(\dfrac{\theta - r}{1 - r}\right), & \theta \in \,]r, 1]. \end{cases} \tag{7.1}$$

Notice that a.s. $U_r(\cdot)$ is nonnegative on $[0, 1]$ and is 0 only at time $r$. Moreover, $U_r(\cdot)$ runs the path of $M$ on $[0, r]$ backwards and then runs the path of $\hat{M}$ on $\,]r, 1]$.

The result is the following:

THEOREM 7.1. *For all $\Phi \in C_b^1(H)$ and $h \in D(A)$*

$$\mathbb{E}[\partial_h \Phi(Y) 1_{(Y \in K)}] = -\mathbb{E}[(\langle Y, h'' \rangle_L - \bar{Y} \cdot \bar{h}) \Phi(Y) 1_{(Y \in K)}] \\ - \int_0^1 h_r \frac{1}{\sqrt{2\pi^3 r(1 - r)}} \mathbb{E}[\Phi(U_r) e^{-(1/2)(\overline{U}_r)^2}] \, dr. \tag{7.2}$$

The proof is postponed to Appendix A and is based on techniques introduced in [24] and [4], where similar results were proven for the law of the Brownian motion and the Brownian bridge conditioned to be greater or equal than a fixed value. See also [16] for related results.

Notice that (7.2) is an integration by parts formula for the law of $Y$ on the set $K$, since the left hand side contains an integration of a partial derivative of $\Phi$, while in the right-hand side only $\Phi$ appears but none of its derivatives. The second term in the right-hand side of (7.2) is an infinite-dimensional boundary term. Indeed, a.s. the typical path of $Y$ conditioned on $K$ is positive on $[0, 1]$. Instead, a.s. $U_r(\cdot)$ is nonnegative and equal to 0 at (and only at) $\theta = r$: therefore it lies in the "boundary" of the set $K$, support of the measure in the left-hand side; see [24].

We denote by $p_{\overline{U}_r}(c)$, $c > 0$, the continuous version of the density of $\overline{U}_r$. By conditioning on $\overline{Y} = c$, we obtain from (7.2):

COROLLARY 7.2. *For all $c > 0$, $h \in D(A)$ and $\Phi \in C_b^1(H)$,*

$$\mathbb{E}[\partial_{\Pi h}\Phi(Y^c)1_{(Y^c \in K)}]$$

(7.3)
$$= -\mathbb{E}[\langle Y^c, h''\rangle_L \ \Phi(Y^c)1_{(Y^c \in K)}]$$

$$- \int_0^1 \Pi h_r \frac{p_{\overline{U}_r}(c)}{\pi\sqrt{r(1-r)}} \ \mathbb{E}[\Phi(U_r)|\overline{U}_r = c] \, dr.$$

Notice that the density of $p_{\overline{U}_r}(c)$ and the law of $U_r$ conditioned on $\{\overline{U}_r = c\}$ are defined in terms of the density of $\langle M, 1\rangle_L$ and the law of $M$ conditioned on $\{\langle M, 1\rangle_L = c\}$; for further information about these objects see [27].

Since $\mu_c$ is the law of $Y^c$ by Lemma 2.1, then we recall that we defined

$$\nu_c(d\omega) = \mathbb{P}(Y^c \in d\omega | Y^c \in K).$$

We also set for all $r \in (0, 1)$, recalling (7.1),

$$\Sigma_r^c(d\omega) := \frac{1}{\mu_c(K)} \frac{p_{\overline{U}_r}(c)}{\pi\sqrt{r(1-r)}} \ \mathbb{P}(U_r \in d\omega | \overline{U}_r = c).$$

Then (7.3) can be rewritten as follows: for all $h \in D(A)$ and $\Phi \in C_b^1(H)$

(7.4) $$\int \partial_{\Pi h}\Phi \, d\nu_c = -\int \langle x, Ah\rangle_L \Phi(x)\nu_c(dx) - \int_0^1 dr \Pi h_r \int \Phi \, d\Sigma_r^c.$$

Notice that we have an analogous formula for $\nu_c^\varepsilon$:

(7.5) $$\int \partial_{\Pi h}\Phi \, d\nu_c^\varepsilon = -\int \langle x, Ah\rangle_L \Phi(x)\nu_c^\varepsilon(dx) - \int_0^1 dr \Pi h_r \int \Phi \, d\Sigma_r^{\varepsilon,c},$$

where

$$\Sigma_r^{\varepsilon,c}(dx) := \frac{1}{\varepsilon} f(x_r)\nu_c^\varepsilon(dx).$$

Then we have that $\int_0^1 dr \Pi h_r \Sigma_r^{\varepsilon,c}$ converges as a measure to $\int_0^1 dr \Pi h_r \Sigma_r^c$ if $\varepsilon \to 0$, as the following lemma states. Notice that this result is crucial in the proof of Proposition 8.1 below, see (8.2).

LEMMA 7.3. *Let $c > 0$.*

1. *For all $\Phi \in C_b(H)$ and $h \in D(A)$,*

$$\lim_{\varepsilon \to 0} \int_0^1 dr \, \Pi h_r \int \Phi \, d\Sigma_r^{\varepsilon,c} = \int_0^1 dr \, \Pi h_r \int \Phi \, d\Sigma_r^c.$$

2. *For all $\delta > 0$ there exists a compact set $C_\delta \subset H_c$ such that $\int_0^1 dr \, \Sigma_r^{\varepsilon,c}(H_c \backslash C_\delta) < \delta$ for all $\varepsilon > 0$.*

The second assertion of this lemma does not follow from the first one and requires a separate proof, because (7.3) and (7.5) contain $\Pi h$ which has zero average, and therefore we cannot compute $\int_0^1 \Sigma_r^{\varepsilon,c}(\cdot) \, dr$ from (7.5).

PROOF OF LEMMA 7.3. Let $\Phi \in C_b^1(H)$. Then the desired convergence holds by (7.5) and dominated convergence, since $1 \geq e^{-U_\varepsilon} \to 1_K$. We recall that $C_b^1(H)$ is dense in $C_b(H)$ in the sup-norm, so that it is enough to prove that

(7.6) $$\limsup_{\varepsilon \to 0} \sup_{\|h\|_\infty \leq 1} \sup_{\|\Phi\|_\infty \leq 1} \left| \int_0^1 dr \, \Pi h_r \int \Phi \, d\Sigma_r^{\varepsilon,c} \right| < \infty.$$

Now, notice that for all $\varepsilon > 0$ and $h \in C([0,1])$,

$$\left| \int_0^1 dr \, \Pi h_r \int \Phi \, d\Sigma_r^{\varepsilon,c} \right| \leq 2 \|h\|_\infty \|\Phi\|_\infty \int_0^1 \Sigma_r^{\varepsilon,c}(1) \, dr.$$

Therefore, (7.6) is proven in particular if we show that for all $\Phi \in C_b(L^2(0,1))$

$$\exists \lim_{\varepsilon \to 0} \int_0^1 \Sigma_r^{\varepsilon,c}(\Phi) \, dr \in \mathbb{R},$$

and this formula is proven in Appendix B. Moreover, notice that this result also yields the second assertion of Lemma 7.3 by Prohorov's theorem. $\square$

**8. The Dirichlet form.** In this section we prove that the process $X$ is associated with a Dirichlet form. The proof is achieved using the integration by parts formula (7.3) and the uniform strong Feller property (3.8) of $X^\varepsilon$. We set for all $\varphi, \psi \in C_b^1(H)$

$$\mathcal{E}^c(\varphi, \psi) := \tfrac{1}{2} \int (\Pi \nabla_H \varphi, \nabla_H \psi)_H \, d\nu_c = \tfrac{1}{2} \int \langle -A \nabla \varphi, \nabla \psi \rangle_L \, d\nu_c.$$

Notice that $\mathcal{E}^c(\varphi, \psi) = \lim_{\varepsilon \to 0} \mathcal{E}^{\varepsilon,c}(\varphi, \psi)$, so it is natural to guess that $\mathcal{E}^c$ is related to equation (1.1). On the other hand, it is not obvious that this is the case: the sole convergence in this sense of a sequence of Dirichlet forms yields essentially no information on the limit, not even that it is also a Dirichlet form. In the next proposition we shall prove first that $\mathcal{E}^c$ is indeed closable and that the resolvent $R^{\varepsilon,c}$ converges as $\varepsilon \to 0$. We refer to [15] and [20] for the general theory of Dirichlet forms.

PROPOSITION 8.1. *Let $c > 0$.*

1. *The term $(\mathcal{E}^c, \mathrm{Exp}_A(H))$ is closable in $L^2(\nu_c)$ and the closure $(\mathcal{E}^c, D(\mathcal{E}^c))$ is a symmetric Dirichlet form such that $\mathrm{Lip}(H_c) \subset D(\mathcal{E}^c)$ and $\mathcal{E}^c(\varphi, \varphi) \leq [\varphi]^2_{\mathrm{Lip}(H_c)}$.*
2. *The term $(P_t^c)_{t \geq 0}$ is the semigroup associated with $(\mathcal{E}^c, D(\mathcal{E}^c))$.*

In particular, $(X_t(x) : t \geq 0, x \in K \cap H_c)$ is associated with the Dirichlet form $\mathcal{E}^c$.

In the proof we use the following easily proven fact.

LEMMA 8.2. *If $E$ is a Polish space, $(M_n : n \in \mathbb{N} \cup \{\infty\})$ is a sequence of finite measures on $E$, such that $\int_E \varphi\, dM_n$ converges to $\int_E \varphi\, dM_\infty$ for all $\varphi \in C_b(E)$, and $(\varphi_n : n \in \mathbb{N} \cup \{\infty\})$ is an equibounded and equicontinuous sequence of functions on $E$, such that $\varphi_n$ converges pointwise to $\varphi_\infty$ on $E$, then $\int_E \varphi_n\, dM_n \to \int_E \varphi_\infty\, dM_\infty$ as $n \to \infty$.*

PROOF OF PROPOSITION 8.1. We divide the proof in several steps.

*Step* 1. We define for all $\varphi \in C_b(H_c)$ the resolvent of $X$ on $K \cap H_c$:

$$R_\lambda^c \varphi(x) := \int_0^\infty e^{-\lambda t} \mathbb{E}[\varphi(X_t(x))]\, dt, \qquad x \in K \cap H_c, \lambda > 0.$$

By (6.1), $R_\lambda^{\varepsilon,c} \varphi(x) \to R_\lambda^c \varphi(x)$ for all $x \in K \cap H_c$. We define $D := \{R_\lambda^c \varphi, \varphi \in C_b(H), \lambda > 0\}$.

Notice first that the integration by parts formula (7.4) can be extended to all $\Phi \in \mathrm{Lip}(H_c)$. Indeed, we can set $\partial_{\Pi h}\Phi := (\Pi \nabla_H \Phi, h)_H$, where $\Pi \nabla_H \Phi \in L^2(H_c, \nu_c; H)$ exists since $\mathrm{Lip}(H_c) \subset D(\mathcal{E}^{\varepsilon,c})$ and $\nu_c$ is absolutely continuous w.r.t. $\mu_c^\varepsilon$. We want to prove now that for all $\psi \in C_b(H)$, the following formula holds:

$$(8.1)\quad \mathcal{E}_\lambda^c(R_\lambda^c \psi, v) := \lambda \int_H R_\lambda^c \psi v\, d\nu_c + \mathcal{E}^c(R_\lambda^c \psi, v) = \int \psi v\, d\nu_c \qquad \forall v \in D.$$

We recall that for $\varphi : H \mapsto \mathbb{C}$, $\varphi(k) := \exp(i(h, k)_H)$, where $h \in D(A^2)$, $i \in \mathbb{C}$ and $i^2 = -1$, the generator $L^\varepsilon$ of the process $X^\varepsilon$, given in (3.9), is

$$L^\varepsilon \varphi(x) := -\frac{1}{2}\varphi(x)\left(\|\Pi h\|_H^2 + i(A^2 h, x)_H - i\frac{\langle f(x), \Pi h\rangle_L}{\varepsilon}\right), \qquad x \in L^2(0, 1).$$

Since $R_\lambda^{\varepsilon,c} = (\lambda - L^\varepsilon)^{-1}$ on $H_c$ is self-adjoint in $L^2(\nu_c^\varepsilon)$,

$$\int R_\lambda^{\varepsilon,c} \psi(\lambda \varphi - L^\varepsilon \varphi)\, d\nu_c^\varepsilon = \int \psi \varphi\, d\nu_c^\varepsilon.$$

Since $e^{-U_\varepsilon} \to 1_K$ as $\varepsilon \to 0$, by Lemma 7.3 we obtain

$$\lim_{\varepsilon \to 0} \int g L^\varepsilon \varphi\, d\nu_c^\varepsilon = \int g L \varphi\, d\nu_c + i \int_0^1 dr\, \Pi h_r \int g\varphi\, d\Sigma_r^c \qquad \forall g \in C_b(H),$$

where $L\varphi$ is defined in (2.8). The crucial fact is now the following: by Lemma 7.3 and Lemma 8.2 we can prove that

$$(8.2) \quad \lim_{\varepsilon \to 0} \int R_\lambda^{\varepsilon,c} \psi(x) \varphi(x) \frac{\langle f(x), \Pi h \rangle_L}{\varepsilon} \nu_c^\varepsilon(dx) = \int_0^1 dr \Pi h_r \int R_\lambda^c \psi \varphi \, d\Sigma_r^c.$$

In particular we obtain again by Lemma 8.2

$$\int \psi \varphi \, d\nu_c = \lim_{i \to \infty} \int \psi \varphi \, d\nu_c^\varepsilon = \lim_{\varepsilon \to 0} \int R_\lambda^{\varepsilon,c} \psi (\lambda \varphi - L^\varepsilon \varphi) \, d\nu_c^\varepsilon$$

$$= \int R_\lambda^c \psi (\lambda \varphi - L\varphi) \, d\nu_c - i \int_0^1 dr \Pi h_r \int R_\lambda^c \psi \varphi \, d\Sigma_r^c$$

and by the integration by parts formula (7.2), the last expression is equal to $\mathcal{E}_\lambda^c(R_\lambda^c \psi, \varphi)$, that is, we have proven (8.1) for $v = \varphi$. By linearity we obtain (8.1) for all $v \in \mathrm{Exp}_A(H)$. By density, we obtain (8.1) for all $v \in D$.

*Step* 2. We want to prove now that the bilinear form $(\mathcal{E}^c, D)$ is closable and the closure is a Dirichlet form. By Lemma I.3.4 in [20], it is enough to prove that if $(u_n)_n \subset D$ and $u_n \to 0$ in $L^2(\nu_c)$ then $\mathcal{E}^c(u_n, v) \to 0$ for any $v \in D$.

By (8.1) we observe that for all $v \in D$ there exists some $\psi_v \in C_b(H)$ such that:

$$\mathcal{E}^c(u, v) = \int u \psi_v \, d\nu_c \qquad \forall \, u \in D.$$

Therefore the above-mentioned criterion applies to $(\mathcal{E}^c, D)$: we denote by $(\tilde{\mathcal{E}}, D(\tilde{\mathcal{E}}))$ the closure. We also obtain that $(R_\lambda^c)_{\lambda > 0}$ is the resolvent of $\tilde{\mathcal{E}}$.

*Step* 3. Finally, we want to show that $(\mathcal{E}^c, \mathrm{Exp}_A(H))$ is closable and that the closure coincides with $(\tilde{\mathcal{E}}, D(\tilde{\mathcal{E}}))$ constructed in the previous step. To this aim it is enough to show that $D(\tilde{\mathcal{E}})$ contains all Lipschitz functions on $K \cap H_c$ and in particular $\mathrm{Exp}_A(H)$. Indeed, the density of $\mathrm{Exp}_A(H)$ follows from the density of this space in $D(\mathcal{E}^{\varepsilon,c})$.

Consider $\psi \in \mathrm{Lip}(H_c) \subset D(\mathcal{E}^{\varepsilon,c})$: by the general theory of Dirichlet forms,

$$(8.3) \qquad \psi \in D(\tilde{\mathcal{E}}) \iff \sup_{\lambda > 0} \int \lambda(\psi - \lambda R_\lambda^c \psi) \psi \, d\nu_c < \infty.$$

By point 1 of Proposition 3.5 we have

$$\int \lambda(\psi - \lambda R_\lambda^{\varepsilon,c} \psi) \psi \, d\nu_c^\varepsilon = \mathcal{E}^{\varepsilon,c}(\lambda R_\lambda^{\varepsilon,c} \psi, \psi) \leq [\psi]_{\mathrm{Lip}(H_c)}^2,$$

so that letting $\varepsilon \to 0$,

$$\int \lambda(\psi - \lambda R_\lambda^c \psi) \psi \, d\nu_c \leq [\psi]_{\mathrm{Lip}(H_c)}^2,$$

and therefore $\mathrm{Lip}(H_c) \subset D(\tilde{\mathcal{E}})$.

In order to prove that $(\mathcal{E}^c, D(\mathcal{E}^c))$ is a Dirichlet form, it is enough to notice that $R^{\varepsilon,c}$ is a Markovian kernel, so that $R^c$ is also Markovian and the result follows from Theorem 4.4 of [20]. $\square$

## APPENDIX A: PROOF OF THEOREM 7.1

**A.1. An absolute continuity result.** To compute an integration by parts formula for $\nu_c$, we want to reduce to a Brownian situation. To this aim, we look for an absolute continuity result between $\mu$ and the law of a Brownian motion with random initial value. We recall the following notation: $B$ is a Brownian motion, $a \in \mathcal{N}(0,1)$ with $\{B,a\}$ independent, and we set $Y_\theta := B_\theta - \overline{B} - a$ for $\theta \in [0,1]$. In Lemma 2.1 we have proven that the law of $Y$ is $\mu$. Then we have the following result:

PROPOSITION A.1. *For all* $\Phi : C([0,1]) \mapsto \mathbb{R}$ *bounded and Borel*

$$(A.1) \qquad \mathbb{E}[\Phi(Y)] = \mathbb{E}\left[\Phi(b+B)\sqrt{\frac{4}{3}}\exp\left(-\frac{1}{2}(b+\overline{B})^2 + \frac{3}{8}b^2\right)\right]$$

*where $B$ is a Brownian motion with $B_0 = 0$, $b \sim \mathcal{N}(0, 4/3)$ and $\{b, B\}$ are independent.*

PROOF. The thesis follows if we show that the Laplace transforms of the two probability measures in (A.1) are equal. Notice first that for all $h \in L^2(0,1)$

$$\mathbb{E}[e^{\langle Y,h \rangle_L}] = e^{1/2 \langle \overline{Q}h, h \rangle_L};$$

recall (2.5), (1.11) and (1.9). Recall now the following version of the Cameron–Martin theorem: for all $\Phi \in C_b(L^2(0,1))$ and $h \in L^2(0,1)$

$$\mathbb{E}[\Phi(B)e^{\langle B,h \rangle_L}] = e^{1/2 \langle Q_B h, h \rangle_L} \mathbb{E}[\Phi(B + Q_B h)], \qquad Q_B h(\theta) := \int_0^1 \theta \wedge \sigma h_\sigma \, d\sigma,$$

and the following standard Gaussian formula for $X \sim N(0, \sigma^2)$, $\sigma \geq 0$ and $\beta \in \mathbb{R}$:

$$\mathbb{E}[e^{-(1/2)(X+\beta)^2}] = \frac{1}{\sqrt{1+\sigma^2}} e^{-(1/2)(\beta^2/(1+\sigma^2))}.$$

Applying these formulae and recalling that $\overline{B} \sim N(0, 1/3)$, we obtain

$$\mathbb{E}\left[e^{\langle b+B, h \rangle_L}\sqrt{\frac{4}{3}}e^{-(1/2)(b+\overline{B})^2 + (3/8)b^2}\right]$$

$$= \frac{1}{\sqrt{2\pi}} \int_\mathbb{R} \mathbb{E}[e^{\langle B,h \rangle_L - (1/2)(y+\overline{B})^2}] e^{y\bar{h}} \, dy$$

$$= \frac{1}{\sqrt{2\pi}} \int_\mathbb{R} \mathbb{E}[e^{-(1/2)(y+\overline{B}+\overline{Q_B h})^2}] e^{y\bar{h} + (1/2)\langle Q_B h, h \rangle_L} \, dy$$

$$= \sqrt{\frac{3}{8\pi}} \int_\mathbb{R} e^{-(3/8)(y+\overline{Q_B h})^2 + y\bar{h} + 1/2 \langle Q_B h, h \rangle_L} \, dy$$

$$= e^{(1/2)(\langle Q_B h, h \rangle_L - 2\bar{h} \cdot \overline{Q_B h} + (4/3)(\bar{h})^2)} = e^{1/2 \langle \overline{Q}h, h \rangle_L}. \qquad \square$$

**A.2. Integration by parts.** We now want to prove Theorem 7.1, where the following integration by parts formula is stated: for all $\Phi \in C_b^1(H)$ and $h \in D(A)$

$$
\begin{aligned}
\mathbb{E}[\partial_h \Phi(Y) 1_{(Y \in K)}] = &-\mathbb{E}[(\langle Y, h'' \rangle_L - \overline{Y} \cdot \bar{h}) \Phi(Y) 1_{(Y \in K)}] \\
&- \int_0^1 h_r \frac{1}{\sqrt{2\pi^3 r(1-r)}} \mathbb{E}[\Phi(U_r) e^{-(1/2)(\overline{U}_r)^2}] \, dr.
\end{aligned}
\tag{A.2}
$$

Denoting by $(M, \hat{M})$ two independent copies of the standard Brownian meander (see [22]), we set for all $r \in (0, 1)$

$$
V_r(\theta) := \begin{cases} -\sqrt{r} M(1) + \sqrt{r} M\left(\dfrac{r-\theta}{r}\right), & \theta \in [0, r], \\ -\sqrt{r} M(1) + \sqrt{1-r} \hat{M}\left(\dfrac{\theta-r}{1-r}\right), & \theta \in \,]r, 1]. \end{cases}
\tag{A.3}
$$

Notice that $V_r(\cdot)$ is 0 at time 0, then runs backward the path of $M$ on $[0, r]$ and then runs the path of $\hat{M}$ on $]r, 1]$. Almost surely since $M > 0$ on $]0, 1]$, then $V_r(\cdot)$ attains the minimum $-\sqrt{r} M(1)$ only at time $\theta = r$. Recalling (7.1), we have

$$
V_r = -\sqrt{r} M(1) + U_r, \qquad r \in [0, 1].
$$

We recall now the following path decomposition of a Brownian motion $B$ on the time interval $[0, 1]$:

THEOREM A.2 ([12]). *Let $(\tau, M, \hat{M})$ be a independent triple, such that $\tau$ has the arcsine law, $M$ and $\hat{M}$ are two standard Brownian meanders. Then $V_\tau \stackrel{d}{=} B$, where $V$ is defined by* (A.3).

PROOF OF THEOREM 7.1. Notice that, integrating out the variable $b \sim N(0, 4/3)$, we can write (A.1) in the following way:

$$
\begin{aligned}
\mathbb{E}[\Phi(Y)] &= \int_\mathbb{R} \mathbb{E}\left[\Phi(y+B) \sqrt{\frac{4}{3}} e^{-(1/2)(y+\overline{B})^2 + (3/8)y^2}\right] \frac{e^{-(3/8)y^2}}{\sqrt{2\pi 4/3}} \, dy \\
&= \int_\mathbb{R} \mathbb{E}_y[\Phi(B) \rho(B)] \, dy,
\end{aligned}
$$

where, under $\mathbb{P}_y$, $B$ is a Brownian motion with $B_0 = y$ and we set $\rho : C([0, 1]) \mapsto \mathbb{R}$,

$$
\rho(\omega) := \frac{1}{\sqrt{2\pi}} \exp\left(-\frac{1}{2}(\bar{\omega})^2\right), \qquad \bar{\omega} := \int_0^1 \omega.
$$

In particular we can write

$$
\mathbb{E}[\partial_h \Phi(Y)] = \int_\mathbb{R} \mathbb{E}_y[\partial_h \Phi(B) \rho(B)] \, dy.
$$

Now, without loss of generality, we can suppose that

$$h \geq 0, \qquad \Phi \geq 0.$$

In particular, $K \subseteq K - th$ for all $t \geq 0$. Recall that $\partial_h \Phi(x) = \lim_{t \downarrow 0}(\Phi(x) - \Phi(x - th))/t$. By the Cameron–Martin theorem

$$\frac{1}{t}\mathbb{E}_y[1_K(B)\rho(B)(\Phi(B) - \Phi(B - th))]$$

$$= -\frac{1}{t}\mathbb{E}_y[1_{(K-th)\setminus K}(B)\Phi(B)\rho(B)]$$

$$- \frac{1}{t}\mathbb{E}_y[1_{K-th}(B)\Phi(B)(\rho(B+th) - \rho(B))]$$

$$+ \frac{1}{t}\mathbb{E}_y\left[1_{K-th}(B)\rho(B+th)\Phi(B)\left(1 - \exp\left(-\frac{1}{2}\|th'\|^2 + t\langle B, h''\rangle\right)\right)\right].$$

Notice that only the limit of the first term in the right-hand side of the last formula is not trivial. By Theorem A.2 we have

$$\mathbb{E}_y[1_{(K-th)\setminus K} \cdot \Phi \cdot \rho(B)] = \mathbb{E}[1_{(K-th)\setminus K} \cdot \Phi \cdot \rho(y + V_\tau)]$$

$$= \int_0^1 \frac{1}{\pi\sqrt{r(1-r)}}\mathbb{E}[1_{(K-th)\setminus K} \cdot \Phi \cdot \rho(y+V_r)]\,dr.$$

Before proceeding, we turn to a similar computation where $h$ is not continuous but a step function, that is, of the form $\sum_{j=1}^n c_j 1_{I_j}$, where $c_j \in \mathbb{R}$ and $I_j \subset [0,1]$ is measurable. Since such functions are dense in $C([0,1])$, this will yield the desired convergence for $h \in C([0,1])$ by density (see the end of the proof).

Let $n \in \mathbb{N}$, $c_n \geq c_{n-1} \geq \cdots \geq c_1 \geq c_0 := 0$, $\{I_1, \ldots, I_n\}$ a Borel partition of $[0,1]$ and $I_0 := \varnothing$, and set

$$h_i := \sum_{j=1}^n (c_j \wedge c_i)1_{I_j}, \qquad i = 0, \ldots, n.$$

The key point is the following: for $i = 1, \ldots, n$, since $h_i \geq h_{i-1}$, and $h_i = h_{i-1}$ on $\bigcup_{j=0}^{i-1} I_j$, then for all $r \in (0,1)$ we have

$$y + V_r \in (K - th_i)\setminus(K - th_{i-1})$$

if and only if

$$V_r \geq -th_i - y \quad \text{on } [0,1], \qquad r \in \bigcup_{j=i}^n I_j \text{ and } \sqrt{r}M(1) \in [y + tc_{i-1}, y + tc_i).$$

Indeed, $V_r$ attains its minimum $-\sqrt{r}M(1)$ only at time $r$. Then we obtain for all $t \geq 0$ and $i = 1, \ldots, n$

$$\mathbb{E}_y[1_{(K-th_i)\backslash K} \cdot \Phi \cdot \rho(B)]$$

$$= \int_0^1 \frac{dr}{\pi\sqrt{r(1-r)}} \mathbb{E}[\Phi \cdot \rho \cdot [1_{(K-th_{i-1})\backslash K} + 1_{(K-th_i)\backslash(K-th_{i-1})}](y+V_r)]$$

$$= \int_0^1 \frac{dr}{\pi\sqrt{r(1-r)}} \mathbb{E}[\Phi \cdot 1_{(K-th_{i-1})\backslash K}(y+V_r)]\, dr$$

$$+ \int_{\bigcup_{j=i}^n I_j} \frac{dr}{\pi\sqrt{r(1-r)}} \mathbb{E}\left[\Phi \cdot 1_{(K-th_i)}(y+V_r) 1_{[c_{i-1},c_i)}\left(\frac{\sqrt{r}M(1)-y}{t}\right)\right].$$

Proceeding by induction on $n$ we obtain

$$\mathbb{E}_y[1_{(K-th_n)\backslash K} \cdot \Phi \cdot \rho(B)]$$

$$= \sum_{i=1}^n \sum_{j=i}^n \int_{I_j} \frac{dr}{\pi\sqrt{r(1-r)}}$$

$$\times \mathbb{E}\left[1_{(K-th_i)} \cdot \Phi \cdot \rho\, (y+V_r) 1_{[c_{i-1},c_i)}\left(\frac{\sqrt{r}M(1)-y}{t}\right)\right]$$

so that, since $\sqrt{r}M(1) + V_r = U_r$ by the definitions (7.1) and (A.3),

$$\lim_{t\downarrow 0} \frac{1}{t} \int_\mathbb{R} \mathbb{E}_y[1_{(K-th_n)\backslash K} \cdot \Phi \cdot \rho(B)]\, dy$$

$$= \sum_{i=1}^n \sum_{j=i}^n \int_{I_j} \frac{c_i - c_{i-1}}{\pi\sqrt{r(1-r)}} \mathbb{E}[1_K \cdot \Phi \cdot \rho(\sqrt{r}M(1)+V_r)]\, dr$$

$$= \int_0^1 \frac{h_n(r)}{\pi\sqrt{r(1-r)}} \mathbb{E}[\Phi(U_r)\rho(U_r)]\, dr.$$

Set now $I_i := h^{-1}([(i-1)/n, i/n))$, $i \in \mathbb{N}$,

$$f_n := \sum_{i=1}^\infty \frac{i-1}{n} 1_{I_i}, \qquad g_n := \sum_{i=1}^\infty \frac{i}{n} 1_{I_i},$$

where both sums are finite, since $h$ is bounded. Then $f_n \leq h \leq g_n$, $f_n$ and $g_n$ converge uniformly on $[0,1]$ to $h$ as $n \to \infty$ and: $K - tf_n \subseteq K - th \subseteq K - tg_n$, $t \geq 0$. Therefore we obtain by comparison

$$\lim_{t\downarrow 0} \frac{1}{t} \mathbb{E}[1_{(K-th)\backslash K} \cdot \Phi \cdot \rho(B)] = \int_0^1 \frac{h(r)}{\pi\sqrt{r(1-r)}} \mathbb{E}[\Phi(U_r)\rho(U_r)]\, dr.$$

Finally, since $\partial_h \rho(\omega) = -\bar{\omega}\bar{h}\rho(\omega)$, we have proven

$$\mathbb{E}[\partial_h \Phi(Y) \cdot 1_K(Y)] = -\mathbb{E}[\Phi(Y)(\langle Y, h''\rangle_L - \overline{Y} \cdot \bar{h}) 1_K(Y)]$$

$$-\int_0^1 \frac{h(r)}{\pi\sqrt{r(1-r)}} \mathbb{E}[\Phi(U_r)\rho(U_r)]\,dr$$

so that (A.2) is proved. □

## APPENDIX B: PROOF OF LEMMA 7.3

The aim of this section is to conclude the proof of Lemma 7.3. In particular, we want to prove that for all $\Phi \in C_b(L^2(0,1))$

(B.1) $$\exists \lim_{\varepsilon\to 0} \int_0^1 \Sigma_r^{\varepsilon,c}(\Phi)\,dr \in \mathbb{R}.$$

By symmetry, it is enough to prove convergence of $\int_0^{1/2} \Sigma_r^{\varepsilon,c}(\Phi)\,dr$. To this aim, we proceed as in the proof of the first assertion of Lemma 7.3, but restricting the analysis to the path space $L^2(0,1/2)$. The advantage is that, on this space, the processes we consider have no more fixed mean, and in fact we can write now integration by parts formulae for $h \in L^2(0,1/2)$ without the constraint of zero mean; actually, we only need to consider $h \equiv 1$ on $[0,1/2]$; see (B.3).

We prove first an integration by parts formula for the law of $Y^c$ on the path space

$$\hat{K} := \{h \in C([0,1]) : h_\theta \geq 0, \forall \theta \in [0,1/2]\}.$$

We also need an absolute continuity result between the law of $Y^c$ and the law of a Brownian motion with a random initial position, in analogy to (A.1). However, since a.s. the trajectories of $Y^c$ have fixed mean $c$ on $[0,1]$, such absolute continuity can hold only if we restrict to an interval like $[0,1/2]$.

LEMMA B.1. *For all $\Psi : C([0,1/2]) \mapsto \mathbb{R}$ bounded and Borel,*

(B.2) $$\mathbb{E}[\Psi(Y^c)] = \sqrt{32}\mathbb{E}[\Psi(b+B)e^{-12(\gamma(b+B)-c)^2+(3/8)b^2}]$$

*where $b \sim N(0,4/3)$ is independent of $B$ and*

$$\gamma(\omega) := \int_0^{1/2}(\omega_r + \omega_{1/2})\,dr, \qquad \omega \in C([0,1/2]).$$

*Moreover, for all $c > 0$ and $\Psi \in C_b^1(L^2(0,1/2))$,*

(B.3) $$\mathbb{E}[\partial_1 \Psi(Y^c)1_{(Y^c \in \hat{K})}] = \mathbb{E}[24(\gamma(Y^c)-c)\Psi(Y^c)1_{(Y^c \in \hat{K})}]$$
$$-\int_0^{1/2}\sqrt{\frac{12}{\pi^3 r(1/2-r)}}\mathbb{E}[\Psi(T_r)e^{-12(\gamma(T_r)-c)^2}]\,dr$$

*where we set for $\theta \in [0,1]$ and $r \in [0, 1/2]$*

$$\mathbf{1}_\theta := \begin{cases} 1, & \theta \in [0, 1/2], \\ 0, & \theta \in \,]1/2, 1], \end{cases}$$

$$T_r(\theta) := \begin{cases} \sqrt{r} M\left(\dfrac{r-\theta}{r}\right), & \theta \in [0, r], \\ \sqrt{\dfrac{1}{2} - r}\, \hat{M}\left(\dfrac{\theta - r}{1/2 - r}\right), & \theta \in \,]r, 1/2]. \end{cases}$$

PROOF. Let us go back to (A.1) and choose $\Phi(\omega) = \Psi(\omega)g(\bar{\omega})$, with $g \in C_b(\mathbb{R})$. Let us notice that

$$b + \overline{B} = \gamma(b+B) + m, \qquad m := \int_{1/2}^{1} (B_r - B_{1/2})\, dr,$$

where $\{(b+B_\theta)_{\theta \in [0,1/2]}, m\}$ are independent and $m \sim N(0, 1/24)$. Then by (A.1)

$$\mathbb{E}[\Psi(Y)g(\overline{Y})]$$
$$= \mathbb{E}[\Phi(Y)]$$
$$= \sqrt{\frac{4}{3}} \mathbb{E}[\Psi(b+B) g(\gamma(b+B)+m) e^{(-1/2)(\gamma(b+B)+m)^2 + (3/8)b^2}]$$
$$= \sqrt{\frac{4}{3}} \int_{\mathbb{R}} \sqrt{\frac{24}{2\pi}} e^{-12y^2} \mathbb{E}[\Psi(b+B) g(\gamma(b+B)+y)$$
$$\qquad \times e^{(-1/2)(\gamma(b+B)+y)^2 + (3/8)b^2}]\, dy$$
$$= \sqrt{32} \int_{\mathbb{R}} g(z) \frac{1}{\sqrt{2\pi}} e^{-1/2 z^2} \mathbb{E}[\Psi(b+B) e^{-12(\gamma(b+B)-z)^2 + (3/8)b^2}]\, dz,$$

and (B.2) is proven. For the proof of (B.3), we write, as in the proof of Theorem 7.1,

$$\mathbb{E}[\Psi(Y^c)] = \int_{\mathbb{R}} \mathbb{E}_y[\Phi(B)\hat{\rho}(B)]\, dy, \qquad \hat{\rho}(\omega) := \sqrt{\frac{12}{\pi}} e^{(-1/2)(\gamma(\omega) - c)^2}.$$

Since $\mathbf{1} \geq 0$, then $\hat{K} \subseteq \hat{K} - t\mathbf{1}$ for all $t \geq 0$ and

$$\frac{1}{t} \mathbb{E}_y[\mathbf{1}_{\hat{K}}(B) \hat{\rho}(B) (\Psi(B) - \Psi(B - t\mathbf{1}))]$$
$$= -\frac{1}{t} \mathbb{E}_y[\mathbf{1}_{(\hat{K} - t\mathbf{1}) \setminus \hat{K}}(B) \Psi(B) \hat{\rho}(B)]$$
$$\quad - \frac{1}{t} \mathbb{E}_y[\mathbf{1}_{\hat{K} - t h}(B) \Psi(B) (\hat{\rho}(B+t) - \hat{\rho}(B))].$$

By Theorem A.2 for a Brownian motion on $[0, 1/2]$ we have

$$\mathbb{E}_y[1_{(\hat{K}-t\mathbf{1})\setminus\hat{K}} \cdot \Psi \cdot \hat{\rho}(B)]$$

$$= \int_0^{1/2} \frac{dr}{\pi\sqrt{r(1/2-r)}} \mathbb{E}[1_{(\hat{K}-t\mathbf{1})\setminus\hat{K}} \Psi \hat{\rho}(y - \sqrt{r}M(1) + T_r)]$$

$$= \int_0^{1/2} \frac{dr}{\pi\sqrt{r(1/2-r)}} \mathbb{E}\left[1_{(\hat{K}-t\mathbf{1})} \Psi \hat{\rho}(y - \sqrt{r}M(1) + T_r) \right.$$

$$\left. \times 1_{[0,1)}\left(\frac{\sqrt{r}M(1) - y}{t}\right)\right]$$

so that

$$\lim_{t\downarrow 0} \frac{1}{t} \int_{\mathbb{R}} \mathbb{E}_y[1_{(\hat{K}-t\mathbf{1})\setminus\hat{K}} \cdot \Psi \cdot \hat{\rho}(B)] \, dy = \int_0^{1/2} \frac{dr}{\pi\sqrt{r(1/2-r)}} \mathbb{E}[\Psi(T_r)\hat{\rho}(T_r)].$$

Finally, since $\partial_1 \hat{\rho}(\omega) = -24(\gamma(\omega) - c)\hat{\rho}(\omega)$, we have proven (A.2). □

We set now for $\varepsilon > 0$, $c > 0$, $r \in (0, 1/2)$ and $\Psi \in C_b(L^2(0, 1/2))$

$$\hat{U}_\varepsilon(x) := \frac{1}{\varepsilon} \int_0^{1/2} F(x(\theta)) \, d\theta, \qquad x \in L^2(0, 1),$$

$$\hat{\Sigma}_r^{\varepsilon,c}(\Psi) := \frac{1}{\varepsilon} \mathbb{E}[f(Y_r^c) \Psi(Y^c) e^{-\hat{U}_\varepsilon(Y^c)}],$$

$$\hat{\Sigma}_r^c(\Psi) := \sqrt{\frac{12}{\pi^3 r(1/2-r)}} \mathbb{E}[\Psi(T_r) e^{-12(\gamma(T_r)-c)^2}].$$

LEMMA B.2. *Let $c > 0$. For all $\Psi \in C_b(L^2(0, 1/2))$,*

$$\lim_{\varepsilon \to 0} \int_0^{1/2} \hat{\Sigma}_r^{\varepsilon,c}(\Psi) \, dr = \int_0^{1/2} \hat{\Sigma}_r^c(\Psi) \, dr.$$

PROOF. Notice that $\nabla \hat{U}_\varepsilon(x) = \frac{1}{\varepsilon} f(x)\mathbf{1}$. By (B.2), for all $\Psi \in C_b^1(L^2(0, 1/2))$,

(B.4)
$$\mathbb{E}[\partial_1 \Psi(Y^c) e^{-\hat{U}_\varepsilon(Y^c)}]$$
$$= \mathbb{E}[24(\gamma(Y^c) - c)\Psi(Y^c) e^{-\hat{U}_\varepsilon(Y^c)}] - \int_0^{1/2} \hat{\Sigma}_r^{\varepsilon,c}(\Psi) \, dr.$$

Comparing this formula with (B.3) we obtain the thesis for all $\Psi \in C_b^1(L^2(0, 1/2))$. Now, by (B.4), for $\Psi \equiv 1$,

$$\int_0^{1/2} \hat{\Sigma}_r^{\varepsilon,c}(1) \, dr = \mathbb{E}[24(\gamma(Y^c) - c) e^{-\hat{U}_\varepsilon(Y^c)}],$$

and therefore
$$\lim_{\varepsilon \to 0} \int_0^{1/2} \hat{\Sigma}_r^{\varepsilon,c}(1) \, dr = \mathbb{E}[24(\gamma(Y^c) - c)1_{(Y^c \in \hat{K})}] < \infty.$$

By density, this yields the thesis for all $\Psi \in C_b(L^2(0, 1/2))$. $\square$

We now set for $\varepsilon > 0$
$$\hat{U}'_\varepsilon(x) := \frac{1}{\varepsilon} \int_{1/2}^1 F(x(\theta)) \, d\theta = U_\varepsilon(x) - \hat{U}_\varepsilon(x), \qquad x \in L^2(0,1).$$

We notice now that, by (A.1) and (B.2), we can compute explicitly the conditional distribution of $Y^c$ given $(Y_r^c, r \in [0, 1/2])$. Indeed, we have for all $\omega \in C([0, 1/2])$ and $\Phi \in C_b(L^2(0,1))$
$$\mathbb{E}[\Phi(Y^c)|Y_\theta^c = \omega_\theta, \forall \theta \in [0, 1/2]] = \mathbb{E}[\Phi(\hat{B}(c, \omega))]$$

where
$$\hat{B}_\theta(c, \omega) := \begin{cases} \omega_\theta, & \theta \in [0, 1/2], \\ \omega_{1/2} + \hat{B}_{\theta-1/2} - \rho_{\theta-1/2}\left(\int_0^{1/2} \hat{B} + \gamma(\omega) - c\right), \\ & \theta \in ]1/2, 1], \end{cases}$$

and $\rho_\theta := 12\theta(1-\theta)$, $\theta \in [0, 1/2]$. Then we have
$$\Sigma_r^{\varepsilon,c}(\Phi) = \frac{1}{Z_c^\varepsilon} \int \mathbb{E}[\Phi(\hat{B}(c,\omega))e^{-\hat{U}'_\varepsilon(\hat{B}(c,\omega))}]\hat{\Sigma}_r^{\varepsilon,c}(d\omega) \qquad \forall r \in [0, 1/2].$$

Since $e^{-\hat{U}'_\varepsilon}$ converges monotonically to $1_{\hat{K}'}$ with $\hat{K}' := \{\omega \in C([0,1]): \omega_\theta \geq 0, \forall \theta \in [1/2, 1]\}$, then it is easy to conclude from Lemma B.2 that $\int_0^{1/2} \Sigma_r^{\varepsilon,c} \, dr$ converges as $\varepsilon \to 0$ to a finite positive measure and therefore it is tight. By symmetry, this yields (B.1).

CONSERVATIVE STOCHASTIC CAHN–HILLIARD EQUATION 35

ENS Cachan
Antenne de Bretagne
Avenue Robert Schumann
35170 Bruz
France
E-mail: [Arnaud.Debussche@bretagne.ens-cachan.fr](mailto:Arnaud.Debussche@bretagne.ens-cachan.fr)

Laboratoire de Probabilités
  et Modèles Aléatoires
Université Paris 6—Pierre et Marie Curie
4 place Jussieu 75252
Paris cedex 05
France
E-mail: [zambotti@ccr.jussieu.fr](mailto:zambotti@ccr.jussieu.fr)
URL: [http://www.proba.jussieu.fr/~zambotti/](http://www.proba.jussieu.fr/~zambotti/)